\documentclass[12pt,reqno]{amsart}%
\usepackage{txfonts}

\usepackage{amsmath, amsfonts, amssymb, amsthm, amscd, amsbsy, latexsym, amsxtra}  \usepackage{fancyhdr}
\usepackage[usenames,dvipsnames,svgnames,x11names,hyperref]{xcolor}
\usepackage{geometry}
\usepackage{graphicx}
\usepackage[pagebackref]{hyperref}
\usepackage{showidx}
\usepackage{amsmath}
\usepackage{amsfonts}
\usepackage{amssymb}%
\providecommand{\U}[1]{\protect\rule{.1in}{.1in}}
\newcommand{\RR}{\mathbb{R}}

\newcommand{\floor}[1]{\lfloor #1 \rfloor}

\hypersetup{colorlinks,breaklinks,
	linkcolor={Fuchsia},
	citecolor={ForestGreen},
	urlcolor={NavyBlue}}
\newtheorem{theorem}{Theorem}[section]
\theoremstyle{plain}

\newtheorem{definition}{Definition}[section]

\newtheorem{lemma}{Lemma}[section]

\newtheorem{proposition}{Proposition}[section]
\newtheorem{remark}{Remark}[section]

\numberwithin{equation}{section}
\allowdisplaybreaks
\def\Xint#1{\mathchoice
	{\XXint\displaystyle\textstyle{#1}}	{\XXint\textstyle\scriptstyle{#1}}	{\XXint\scriptstyle\scriptscriptstyle{#1}}	{\XXint\scriptscriptstyle\scriptscriptstyle{#1}}	\!\int}
\def\XXint#1#2#3{{\setbox0=\hbox{$#1{#2#3}{\int}$ }
		\vcenter{\hbox{$#2#3$ }}\kern-.6\wd0}}
\def\ddashint{\Xint=}
\def\dashint{\Xint-}
\begin{document}
	
	\title[]{Gagliardo--Nirenberg type inequalities using Fractional Sobolev spaces and Besov spaces}
	\author{Nguyen Anh Dao}
	\address{Nguyen Anh Dao: School of Economic Mathematics and Statistics, University of Economics Ho Chi Minh City (UEH),  Vietnam}
	\email{anhdn@ueh.edu.vn}
	
	\date{\today}

	\begin{abstract} Our main purpose is to establish Gagliardo--Nirenberg type inequalities using fractional homogeneous Sobolev spaces, and homogeneous Besov spaces. In particular, we extend some of the results obtained by the authors in \cite{Brezis1, Brezis2, Brezis3,  DaoLamLu1,  Miyazaki, Van}.
	\end{abstract}
	\maketitle

	\section{Introduction}
	In this paper, we are interested in
	the following Gagliardo--Nirenberg inequality:
	\\
	For every $0\leq  \alpha_1<\alpha_2$, and for  $1\leq p_1, p_2, q \leq  \infty$, there holds 
	\begin{equation}\label{-10}
		\|f\|_{\dot{W}^{\alpha_1, p_1}}  \lesssim \| f \|^{1-\frac{\alpha_1}{\alpha_2}}_{ L^{q} }   \|f \|^\frac{\alpha_1}{\alpha_2}_{\dot{W}^{\alpha_2, p_2}}
		\,,
	\end{equation}
	where $$\frac{1}{p_1} =  \frac{1}{q} \left(1-\frac{\alpha_1}{\alpha_2}\right) + \frac{1}{p_2} \frac{\alpha_1}{\alpha_2}  \,,$$ 
	and $\dot{W}^{\alpha,p}(\RR^n)$ denotes by the homogeneous Sobolev space (see its definition in Section 2).
	\\
	It is known that such an inequality of this type  plays an important role in the analysis of PDEs.  When $\alpha_i$, $i=1,2$ are  nonnegative integer  numbers, \eqref{-10} was obtained independently  by  Gagliardo  \cite{Gag} and Nirenberg \cite{Nir}.
	After that,  the inequalities of this type   have been studied by many authors in \cite{Brezis1, Brezis2, Brezis3, CDDD,  Dao1,DaoLamLu1,DaoLamLu2,Le, LuWheeden,Lu2, Miyazaki,MeRi2003, Van}, and the references cited  therein.  
	\\  
	The case $q=\infty$ can be considered as a limiting case of \eqref{-10}, i.e:
	\begin{equation}\label{-12}
		\|D^{\alpha_1}  f\|_{L^{p_1}}  \lesssim  	\|f\|^{1-\frac{\alpha_1}{\alpha_2}} _{L^\infty}  \|D^{\alpha_2}  f  \|^{\frac{\alpha_1}{\alpha_2}}_{L^{p_2}}  \,,  \quad \forall  f\in L^\infty(\RR^n)\cap \dot{W}^{\alpha_2,p_2}(\RR^n)\,,
	\end{equation}
	with $p_1=\frac{p_2\alpha_2}{\alpha_1}$.
	Obviously, this inequality fails if $\alpha_1=0$. \\
	An partial  improvement of \eqref{-12} in terms of {\rm BMO} space  was obtained by Meyer and Rivi\`ere \cite{MeRi2003}:
	\begin{equation}\label{-13}
		\|D f\|^2_{L^{4}}  \lesssim   \|f\|_{ \rm{BMO} } \| D^2 f\|_{ L^2 } \,,
	\end{equation}
	for all $f\in {\rm BMO}(\mathbb{R}^n) \cap W^{2,2}(\mathbb{R}^n)$.  Thanks to \eqref{-13}, the authors proved a regularity result for a class of stationary Yang--Mills fields in high dimension.
	\\
	After that \eqref{-14} was extended to higher derivatives by the authors in \cite{Strz,Miyazaki}.  Precisely, there holds true
	\begin{equation}\label{-14}
		\|D^{\alpha_1}  f\|_{L^{p_1}}  \lesssim  	\|f\|^{1-\frac{\alpha_1}{\alpha_2}} _{{\rm BMO}}  \|D^{\alpha_2}  f  \|^{\frac{\alpha_1}{\alpha_2}}_{L^{p_2}}  \,,
	\end{equation}
	for all $f\in {\rm BMO}(\RR^n) \cap    W^{\alpha_2,p_2}(\mathbb{R}^n)$, $p_2>1$.
	\\
	Recently, the author et al. \cite{DaoLamLu1} improved \eqref{-14} by means of the homogeneous Besov spaces. For convenience, we recall the result here.
	\begin{theorem}[see Theorem 1.2,  \cite{DaoLamLu1}] \label{TheC}     \sl
		Let  $m, k$ be integers with $1\leq  k<m$. For every   $s\geq 0$,  let $f \in \mathcal{S}'(\mathbb{R}^n)$ be such that
		$ D^m  f\in  L^{p}(\mathbb{R}^n)$,  $1\leq p<\infty$; and $f\in\dot{B}^{-s}(\mathbb{R}^n)$. Then,   we have  $D^k  f\in L^r(\mathbb{R}^n)$, $r=p \left( \frac{m+s}{k+s}  \right)$, and
		\begin{equation}\label{-15}
			\|D^k  f\|_{L^r}   \lesssim  \|f\|^{\frac{m-k}{m+s}}_{\dot{B}^{-s}}
			\left\|D^m f \right\|^\frac{k+s}{m+s}_{L^p} \,,
		\end{equation}
		where we denote  $\dot{B}^{\sigma}  =  \dot{B}^{\sigma,\infty}_{\infty}$, $\sigma\in\mathbb{R}$    (see the definition of Besov spaces in Section 2).	
	\end{theorem}
	\begin{remark}  
		Obviously, \eqref{-15} is stronger than \eqref{-14} when $s=0$ since ${\rm BMO}(\mathbb{R}^n)  \hookrightarrow  \dot{B}^{0}(\mathbb{R}^n)$. We emphasize that \eqref{-15} is still true for $k=0$ when $s>0$.
	\end{remark}
	\begin{remark}
		In studying the space ${\rm BV}(\RR^2)$, A. Cohen et al., \cite{CDPX} 
		proved \eqref{-15} for the case $k=0, m=p=1, s=n-1, r=\frac{n}{n-1}$ by using wavelet decompositions (see  \cite{Le} for the case  $k=0, m=1, p\geq 1,  r=p\big(\frac{1+s}{s}\big)$, with $s>0$). 
	\end{remark}
	Inequality \eqref{-10} in terms of fractional Sobolev spaces  has been investigated by the authors in \cite{Brezis1, Brezis2, Brezis3,Van} and the references therein. Surprisingly, there is a border line for the limiting case of Gagliardo--Nirenberg type inequality. In \cite{Brezis1}, Brezis--Mironescu proved that  the following inequality  
	\begin{align}\label{-16}
		\|f\|_{W^{\alpha_1,p_1}}\lesssim 	\|f\|^\theta_{W^{\alpha,p}}	\|f\|^{1-\theta}_{W^{\alpha_2,p_2}}    \,,
	\end{align}
	with  $\alpha_1=\theta \alpha +(1-\theta)\alpha_2$, $\frac{1}{p_1}=\frac{\theta}{p}+\frac{1-\theta}{p_2}$, and $\theta\in(0,1)$ 
	holds  if and only if   
	\begin{equation}\label{special-cond} \alpha-\frac{1}{p}< \alpha_2-\frac{1}{p_2}  \,.\end{equation}
	As a consequence of this result,
	the following inequality 
	\[
	\|f\|_{\dot{W}^{\alpha_1,p_1}}\lesssim 	\|f\|_{L^\infty}	\|D f\|_{L^1}   
	\]	
	fails whenever $0<\alpha_1<1$. 
	\\
	We note that the limiting case of \eqref{-16}  reads as:
	\begin{equation}\label{-17}
		\|f\|_{\dot{W}^{\alpha_1,p_1}}\lesssim 	\|f\|_{L^\infty}\|f\|_{\dot{W}^{\alpha_2,p_2}}  \,,
	\end{equation}
	where $\alpha_1<\alpha_2$, and $\alpha_1  p_1=\alpha_2 p_2$.
	\\
	When $\alpha_2<1$, Brezis--Mironescu improved  \eqref{-17} by means of  ${\rm BMO}(\RR^n)$ using the Littlewood--Paley decomposition. Very recently, Van Schaftingen \cite{Van}  studied \eqref{-17} for the case $\alpha_2=1$  on a convex open set $\Omega\subset \RR^n$ satisfying certain condition.  Particularly, he proved that 
	\begin{equation}\label{-20}
		\|f\|_{\dot{W}^{\alpha_1,p_1}}\lesssim 	\|f\|^{1-\alpha_1}_{{\rm BMO}}	
		\|D f\|^{\alpha_1}_{L^{p_2}}   
	\end{equation}
	where $0<\alpha_1<1$,  $p_1\alpha_1=p_2$, $p_2>1$.
	\\
	Inspired by the above results, we would like to  study \eqref{-10} by means of fractional Sobolev spaces and Besov spaces. Moreover, we also improve the limiting cases \eqref{-17}, \eqref{-20}  in terms of  $\dot{B}^0(\RR^n)$. 
	\subsection*{Main results}	
	Our first result is to improve \eqref{-10} by using fractional Sobolev spaces, and homogeneous Besov spaces.
	\begin{theorem}\label{Mainthe} Let $\sigma>0$, and  $0\leq \alpha_1<\alpha_2<\infty$. Let  $1\leq p_1, p_2  \leq \infty$ be such that $p_1=p_2\big(\frac{\alpha_2+\sigma}{\alpha_1+\sigma}\big)$, and $p_2(\alpha_2+\sigma)>1$. 	If $f\in \dot{B}^{-\sigma}(\RR^n) \cap \dot{W}^{\alpha_2,p_2}(\RR^n)$, then  $f\in \dot{W}^{\alpha_1,p_1}(\RR^n)$. Moreover,  there is a positive  constant $C=C(n,\alpha_1,\alpha_2,p_2, \sigma)$ such that
		\begin{equation}\label{-3}
			\|f\|_{\dot{W}^{\alpha_1,p_1}}  \leq C 
			\|f\|^{\frac{\alpha_2-\alpha_1}{\alpha_2+\sigma}}_{\dot{B}^{-\sigma}}  \|f\|^{\frac{\alpha_1+\sigma}{\alpha_2+\sigma}}_{\dot{W}^{\alpha_2,p_2}}
			\,.
		\end{equation}
	\end{theorem}
	\begin{remark}  Note that 
		\eqref{-3} is not true for the limiting case  $\sigma= \alpha_1=0,  p_1=\infty$, even  \eqref{special-cond} holds, i.e: $\alpha_2-\frac{1}{p_2}>0$. Indeed, if it is  the case,  then \eqref{-3} becomes
		\[
		\|f\|_{L^{\infty}}  \lesssim \|f\|_{\dot{B}^{0}}  \,.
		\]
		Obviously, 
		the inequality cannot happen  since
		$L^\infty(\RR^n)\hookrightarrow  {\rm BMO}(\RR^n)\hookrightarrow \dot{B}^0(\RR^n)$.
	\end{remark}
	However, if $\alpha_1$ is positive, then
	\eqref{-3}  holds true with $\sigma=0$.  This assertion is in the following theorem.
	\begin{theorem}\label{Mainthe1}  Let 
		$\alpha_2> \alpha_1>0$, and let $1\leq p_1, p_2\leq \infty$ be such that $p_1=\frac{\alpha_2 p_2}{\alpha_1}$, and $\alpha_2 p_2>1$. If $f\in \dot{B}^{0}(\RR^n) \cap \dot{W}^{\alpha_2,p_2}(\RR^n)$, then  $f\in \dot{W}^{\alpha_1,p_1}(\RR^n)$. Moreover,  we have 
		\begin{equation}\label{4.1}
			\|f\|_{\dot{W}^{\alpha_1,p_1}}  \lesssim
			\|f\|^{\frac{\alpha_2-\alpha_1}{\alpha_2}}_{\dot{B}^{0}}  \|f\|^{\frac{\alpha_1}{\alpha_2}}_{\dot{W}^{\alpha_2,p_2}}
			\,.
		\end{equation}
	\end{theorem}
	Our paper is organized as follows. In the next section, we provide the definitions of fractional Sobolev spaces and homogeneous Besov spaces. Section 3 is devoted to the proofs of Theorems \ref{Mainthe}, \ref{Mainthe1}. Moreover, we also obtain the homogeneous version of \eqref{-16} with an elementary proof, see Lemma \ref{Lem-Hom-Sobolev}.
	Finally, we prove $\|f\|_{\dot{W}^{s,p}}\approx \|f\|_{\dot{B}^{s,p}}$ for $0<s<1$, $1\leq p<\infty$  in the Appendix section.

	\section{Definitions and preliminary results}
	\subsection{Fractional Sobolev spaces}
	\begin{definition}\label{Def-frac-Sob}   For any $0<\alpha<1$, and for $1\leq p<\infty$,
		we denote  $\dot{W}^{\alpha,p}(\RR^n)$ (resp. $W^{\alpha,p}(\RR^n)$)  by the homogeneous fractional
		Sobolev space  (resp. the  inhomogeneous fractional Sobolev space), endowed by the semi-norm:
		\[  \|f\|_{\dot{W}^{\alpha,p}}  =
		\left(\displaystyle \int_{\RR^n} \int_{\RR^n} \frac{|f(x+h)-f(x)|^p}{|h|^{n+\alpha  p}} dhdx \right)^{\frac{1}{p}} \,,
		\]	
		and the norm
		\[
		\|f\|_{W^{\alpha,p}}  =  \left(\|f\|^p_{L^p} +  \|f\|^p_{\dot{W}^{\alpha,p}} \right)^\frac{1}{p}\,. 
		\]
	\end{definition}
	When $\alpha\geq 1$,   we can define the  higher order fractional Sobolev space as follows:
	\\
	Denote $\floor{\alpha}$ by the integer part of  $\alpha$. Then, we define
	\[ 
	\|f\|_{\dot{W}^{\alpha,p}} =\left\{	\begin{array}{cl}
		&\|D^{\floor{\alpha}} f\|_{L^p} ,\quad \text{if }\, \alpha\in\mathbb{Z}^+.
		\vspace*{0.1in}\\
		& \|D^{\floor{\alpha}}f\|^p_{\dot{W}^{\alpha-\floor{\alpha},p}}   ,\quad\text{otherwise}\,.
	\end{array}\right.
	\]
	In addition, we also define
	\[ 
	\|f\|_{W^{\alpha,p}} =\left\{	\begin{array}{cl}
		&\|f\|_{W^{\alpha,p}} ,\quad \text{if }\, \alpha\in\mathbb{Z}^+.
		\\
		& \left( \|f\|^p_{W^{\floor{\alpha},p}} +  \|D^{\floor{\alpha}} f\|^p_{\dot{W}^{\alpha-\floor{\alpha},p}}  \right)^{\frac{1}{p}} ,\quad\text{otherwise}\,.
	\end{array}\right.
	\]
	\subsection*{Notation} 	Through the paper, we accept the notation
	$\dot{W}^{\alpha,\infty}(\RR^n)=\dot{C}^{\alpha}(\RR^n)$, $\alpha\in(0,1)$; and  $\dot{W}^{0,p}(\RR^n)=L^p(\RR^n)$, $1\leq p\leq \infty.$
	\\
	In addition, we always denote constant by C, which may change
	from line to line. Moreover, the notation $C(\alpha, p,n)$ means that $C$ merely depends
	on $\alpha, p,n$.
	Next, we write $A \lesssim B$ if there exists a constant $c > 0$ such
	that $A <cB$. And we write $A \approx B$ iff $A \lesssim B \lesssim A$.

	\subsection{Besov spaces}
	To define  the homogeneous Besov spaces, we recall the Littlewood--Paley decomposition (see \cite{Triebel}). Let $\phi_j(x)$  be the inverse Fourier transform of the $j$-th component of the dyadic decomposition i.e.,
	$$\sum_{j\in \mathbb{Z}}  \hat{\phi}(2^{-j} \xi )  =1$$
	except $\xi=0$, where
	$ {\rm supp}( \hat{\phi})\subset \left\{ \frac{1}{2} < |\xi| < 2  \right\}$.
	\\
	Next, let us put
	$$ \mathcal{Z}(\mathbb{R}^n) = \left\{ f \in  \mathcal{S}(\mathbb{R}^n), D^\alpha \hat{f}(0) = 0,\, \forall\alpha \in \mathbb{N}^n, \text{ multi-index}  \right\}  \,,$$
	where  $\mathcal{S}(\mathbb{R}^n)$ is the Schwartz space as usual.
	\begin{definition}\label{Def1} For every $s\in\mathbb{R}$, and for every $1\leq p, q\leq \infty$,  the homogeneous Besov space is denoted by
		$$\dot{B}^s_{p,q} =\left\{ f\in \mathcal{Z}'(\mathbb{R}^n):  \|f\|_{\dot{B}^s_{p,q}}  <\infty   \right\} \,,$$
		with
		$$
		\|f\|_{\dot{B}^s_{p,q}}  =  \left\{  \begin{array}{cl}
			&\left(  \displaystyle \sum_{j\in\mathbb{Z}}
			2^{jsq} \|\phi_j * f\|^q_{L^p} \right)^\frac{1}{q}\,, \text{ if   }\,  1\leq q<\infty,
			\vspace{0.1in}  \\
			&  \displaystyle\sup_{ j \in\mathbb{Z} } \left\{ 2^{js} \|\phi_j * f\|_{L^p} \right\}  \,, \text{ if   }\,  q=\infty \,.
		\end{array} \right.  $$
		When $p=q=\infty$, we denote $\dot{B}^s_{\infty,\infty}=\dot{B}^s$ for short. 
	\end{definition}

	The following characterization of  $\dot{B}^{s}_{\infty,\infty}$ is useful for our proof below.
	\begin{theorem}[see Theorem 4, p. 164, \cite{Peetre}]\label{ThePeetre}  Let $\big\{\varphi_\varepsilon\big\}_\varepsilon$ be a sequence of functions such that
		\[\left\{   
		\begin{array}{cl}
			&{\rm supp}(\varphi_\varepsilon)\subset  B(0,\varepsilon)  ,  \quad \big\{ \frac{1}{2\varepsilon}\leq |\xi|\leq \frac{2}{\varepsilon} \big\}\subset \big\{\widehat{\phi_\varepsilon}(\xi) \not=0  \big\}  ,
			\vspace*{0.1in}\\
			&\int_{\RR^n}  x^\gamma  \varphi_\varepsilon (x)\, dx =0 ,\,  \text{for all  multi-indexes }\,  |\gamma|<k,   \text{ where $k$ is a given integer},
			\vspace*{0.1in}\\
			&  \big|D^\gamma  \varphi_\varepsilon(x)\big|  \leq C \varepsilon^{-(n+|\gamma|)}\,  \text{ for every  multi-index } \gamma\,. 
		\end{array}
		\right.\] 
		Assume $s<k$. Then, we have  
		\[
		f\in \dot{B}^s(\RR^n)  \Leftrightarrow
		\sup_{\varepsilon>0}
		\left\{\varepsilon^{-s} \|\varphi_\varepsilon * f\|_{L^\infty} \right\}   < \infty  \,.
		\]
	\end{theorem}
	We end this section by recall the following result (see \cite{DaoLamLu1}). 
	\begin{proposition}[Lifting operator]\label{Pro1}
		Let $s\in\mathbb{R}$, and let $\gamma$ be a multi-index. Then, 
		$\partial^\gamma$ maps $\dot{B}^s(\RR^n)  \rightarrow  \dot{B}^{s-|\gamma|}(\RR^n)$.
	\end{proposition}

	\section{Proof of the Theorems}
	
	\subsection{Proof of Theorem \ref{Mainthe}} 
	We first prove Theorem \ref{Mainthe}  for the case $0\leq \alpha_1<\alpha_2\leq 1$. After that, we consider $\alpha_i\geq 1$, $i=1,2$.
	\\
	
	{\bf i) Step 1: $0\leq \alpha_1<\alpha_2  \leq  1$.}  We divide our argument into the  following cases.
	\\
	
	{\bf a) The case $p_1=p_2=\infty$,  $0< \alpha_1<\alpha_2 <1$.}  Then, \eqref{-3} becomes
	\begin{equation}\label{1.0}
		\|f\|_{\dot{C}^{\alpha_1}}  \lesssim 	\|f\|^{\frac{\alpha_2-\alpha_1}{\alpha_2+\sigma}}_{\dot{B}^{-\sigma}}  \|f\|^{\frac{\alpha_1+\sigma}{\alpha_2+\sigma}}_{\dot{C}^{\alpha_2}}  \,.
	\end{equation}
	To prove \eqref{1.0}, we use a characterization of homogeneous Besov space   $\dot{B}^{s}$ in Theorem \ref{ThePeetre}, and the fact that  $\dot{B}^s(\RR^n)$ coincides with $\dot{C}^s(\RR^n)$,  $s\in(0,1)$ (see \cite{Grevholm}).  
	\\
	Then, let us recall  sequence  $\{\varphi_\varepsilon\}_{\varepsilon>0}$  in Theorem \ref{ThePeetre}. 	
	\\  
	For $\delta>0$, we write 
	\begin{align}\label{2.-1}
		\varepsilon^{-\alpha_1}\|\varphi_\varepsilon  * f\|_{L^\infty}		 &=\varepsilon^{\alpha_2-\alpha_1}  \varepsilon^{-\alpha_2}\|\varphi_\varepsilon  * f\|_{L^\infty}		\mathbf{1}_{\big\{\varepsilon<\delta\big\}}+ \varepsilon^{-(\alpha_1+\sigma)}  \varepsilon^{\sigma}\|\varphi_\varepsilon  * f\|_{L^\infty}		\mathbf{1}_{\big\{\varepsilon\geq \delta\big\}}  
		\\
		&\leq \delta^{\alpha_2-\alpha_1} \|f\|_{\dot{B}^{\alpha_2}}
		+\delta^{-(\alpha_1+\sigma)} \|f\|_{\dot{B}^{-\sigma}}  \,.  \nonumber
	\end{align}
	Minimizing the  right hand side with respect to $\delta$  in the indicated inequality yields
	\begin{align*}
		\varepsilon^{-\alpha_1}\|\varphi_\varepsilon  * f\|_{L^\infty}\lesssim 	\|f\|_{\dot{B}^{-\sigma}}^{\frac{\alpha_2-\alpha_1}{\alpha_2+\sigma}}
		\|f\|_{\dot{B}^{\alpha_2}}^{\frac{\alpha_1+\sigma}{\alpha_2+\sigma}}\,.		
	\end{align*}
	Since the last inequality holds for every $\varepsilon>0$, then we obtain \eqref{1.0}.	
	\begin{remark}\label{Rem2} It is not difficult to observe that the above proof also adapts to the  two following  cases: 
		\begin{enumerate}
			\item[$\bullet$]    $\alpha_1=0$, $\alpha_2<1$, $\sigma>0$. Then, we have
			\begin{equation}\label{2.-3}
				\|f\|_{L^\infty}	\lesssim 	\|f\|^{\frac{\alpha_2}{\alpha_2+\sigma}}_{\dot{B}^{-\sigma}}
				\|f\|_{\dot{B}^{\alpha_2}}^{\frac{\sigma}{\alpha_2+\sigma}}\,.
			\end{equation}
			\item[$\bullet $] $\alpha_1=0$, $\alpha_2<1$,  $\sigma>0$. Then,  we have
			\begin{equation}\label{2.-2}
				\|f\|_{\dot{B}^{\alpha_1}}	\lesssim 	\|f\|^{\frac{\alpha_2-\alpha_1}{\alpha_2}}_{\dot{B}^{0}}
				\|f\|_{\dot{B}^{\alpha_2}}^{\frac{\alpha_1}{\alpha_2}}\,.
			\end{equation}
		\end{enumerate}		
		This is Theorem \ref{Mainthe1} when $p_i=\infty$, $i=1, 2$.	
	\end{remark}
	To end part {\bf a)}, it remains to prove \eqref{-3} for the case $\alpha_2=1$. That is
	\begin{equation}\label{2.-4}
		\|f\|_{\dot{B}^{\alpha_1}}	\lesssim 	\|f\|^{\frac{1}{1+\sigma}}_{\dot{B}^{-\sigma}}
		\|Df\|_{L^\infty}^{\frac{\sigma}{1+\sigma}}\,.
	\end{equation}
	The proof is similar to the one of \eqref{1.0}. Hence, it suffices to prove that
	\begin{equation}\label{2.-5}
		\varepsilon^{-\alpha_1} \|\varphi_\varepsilon  * f\|_{L^\infty}		\mathbf{1}_{\big\{\varepsilon<\delta\big\}}  \leq   \delta^{1-\alpha_1} \|Df\|_{L^\infty}  \,. \end{equation}
	Indeed, using vanishing moment of $\varphi_\varepsilon$  and the mean value theorem  yields
	\begin{align*}
		\big|\varphi_\varepsilon  * f(x)\big|  &=\big| \int_{B(0,\varepsilon)} (f(x)-f(x-y)) \varphi_\varepsilon (y)\, dy   \big| 
		\\
		&\leq   \int_{B(0,\varepsilon)} \|Df\|_{L^\infty}  |y| |\varphi_\varepsilon (y)|\, dy  \leq \varepsilon  \|\varphi_\varepsilon\|_{L^1} \|Df\|_{L^\infty}  \lesssim    \varepsilon   \|Df\|_{L^\infty} \,.
	\end{align*}
	Thus, \eqref{2.-5} follows easily.
	\\
	By repeating the proof of \eqref{2.-1}, we obtain  \eqref{2.-4}.
	\\

	{\bf b) The case $p_i<\infty, \,i=1,2$.}  Then, the proof follows by way of the following lemmas.
	\begin{lemma}\label{Lem10}
		Let $0<\alpha< 1$, and  $1\leq p<\infty$. 	For every $s>0$, if  $f\in \dot{B}^{-s}(\mathbb{R}^n)\cap \dot{W}^{\alpha,p}(\RR^n)$,  then there exists a positive constant $C=C(s,\alpha,p)$ such that 
		\begin{equation}\label{1.1}
			|f(x)| \leq  C
			\| f \|_{\dot{B}^{-s}}^\frac{\alpha}{s+\alpha}   \big[
			\mathbf{G}_{\alpha,p}(f)(x)\big]^{\frac{s}{s+\alpha}}\,, \quad  \text{for    }  x\in\mathbb{R}^n\,,  
		\end{equation}
		with $$\mathbf{G}_{\alpha,p}(f)(x)= \displaystyle\sup_{\varepsilon>0} \left(\fint_{B(0,\varepsilon)}  \frac{|f(x)-f(x-y)|^p}{\varepsilon^{\alpha p}}  dy  \right)^\frac{1}{p}  \,.$$	
	\end{lemma}	
	\begin{remark}\label{Rem6} When $\alpha=1$, then \eqref{1.1} becomes
		\begin{equation}\label{1.1b}
			|f(x)|\leq C\|f\|_{\dot{B}^{-s}}^\frac{1}{s+1}   \big[\mathbf{M}(|Df|)(x)\big]^{\frac{s}{s+1}}\,, \quad  \text{for    }  x\in\mathbb{R}^n\,.
		\end{equation}
		This inequality was obtained by the authors in \cite{DaoLamLu1}. As a result, we get
		\begin{equation}\label{1.1a}
			\|f\|_{L^{p_1}}  \lesssim 	\| f \|_{\dot{B}^{-s}}^\frac{1}{s+1} 
			\|Df\|_{L^{p_2}}\,,   
		\end{equation}
		with $p_1=p_2\big(\frac{s+1}{s}\big)$, $p_2\geq 1$. 
		\\
		This is also  Theorem \ref{Mainthe} when $\alpha_1=0$,  $\alpha_2=1$, $s=\sigma>0$.
	\end{remark}
	\begin{remark}\label{Rem4}  Obviously, for $1\leq p<\infty$ we have  $\|\mathbf{G}_{\alpha,p}(f)\|_{L^p}\lesssim \|f\|_{\dot{W}^{\alpha,p}}$, and $\mathbf{G}_{\alpha,1}(f)(x)\leq \mathbf{G}_{\alpha,p}(f)(x)$ for  $x\in\RR^n$.  
		\\
		Next, applying Lemma \ref{Lem10}  to  $s=\sigma,   \alpha=\alpha_2$, $p=p_2$,  and taking the $L^{p_1}$-norm of \eqref{1.1} yield
		\[ 	
		\|f\|_{L^{p_1}} \lesssim 
		\|f\|_{\dot{B}^{-\sigma}}^\frac{\alpha_2}{\sigma+\alpha_2} \left( \int
		\big|\mathbf{G}_{\alpha_2,p_2}(f)(x)\big|^{\frac{\sigma p_1}{\sigma+\alpha_2}}  \, dx  \right)^{1/p_1}  \leq \|f\|_{\dot{B}^{-\sigma}}^\frac{\alpha_2}{\sigma+\alpha_2}   \|f\|^{\frac{\sigma}{\sigma+\alpha_2}}_{\dot{W}^{\alpha_2,p_2}} \,,
		\] 
		with $p_1= p_2  \big(\frac{\sigma+\alpha_2}{\sigma}\big)$.
		\\
		Hence, we obtain Theorem \ref{Mainthe} for the case $\alpha_1=0$.
	\end{remark}
	
	
	\begin{proof}[Proof of Lemma \ref{Lem10}] Let us recall  sequence $\{\varphi_\varepsilon\}_{\varepsilon>0}$ above. Then,  we have from the triangle inequality that
		\begin{align*}
			|f(x)| 
			\leq |  \varphi_\varepsilon  * f(x)| + | f(x)-\varphi_\varepsilon  * f(x)| :=\mathbf{I}_1+ \mathbf{I}_2 \,.
		\end{align*} 
		We first estimate  $\mathbf{I}_1$ in terms of $\dot{B}^{-s}$. Thanks to Theorem \ref{ThePeetre},   we get
		\begin{align}\label{1.2}
			\mathbf{I}_1 =	\varepsilon^{-s} \varepsilon^{s} | \varphi_\varepsilon  * f(x)|  \leq C \varepsilon^{-s}  \| f \|_{\dot{B}^{-s}}  \,.
		\end{align}
		For $\mathbf{I}_2$, applying H\"older's inequality yields
		\begin{align}\label{1.3}
			\mathbf{I}_2  &\leq  \int_{B(0,\varepsilon)} |f(x)-f(x-y)| \varphi_\varepsilon (y) \, dy  = \varepsilon^{\frac{n}{p}+\alpha} \int_{B(0,\varepsilon)}  \frac{|f(x)-f(x-y)|}{\varepsilon^{\frac{n}{p}+\alpha}} \varphi_\varepsilon (y) \, dy  \nonumber
			\\
			&\leq \varepsilon^{\frac{n}{p}+\alpha}   \|\varphi_\varepsilon\|_{L^{p'}} \left(\int_{B(0,\varepsilon)}  \frac{|f(x)-f(x-y)|^p}{\varepsilon^{n+\alpha p}}  dy  \right)^\frac{1}{p}  \nonumber
			\\
			&\lesssim \varepsilon^{\frac{n}{p}+\alpha}   \|\varphi_\varepsilon\|_{L^{\infty}} \big|B(0,\varepsilon)\big|^{\frac{1}{p'}}  \mathbf{G}_{\alpha,p}(f)(x) \lesssim \varepsilon^{\alpha}   \mathbf{G}_{\alpha,p}(f)(x)  \,.
		\end{align}
		Note that the last inequality follows by  using the fact $\|\varphi_\varepsilon\|_{L^{\infty}}\leq C \varepsilon^{-n}$.
		\\
		By combining \eqref{1.2} and \eqref{1.3}, we obtain
		\[|f(x)|\leq C\left(\varepsilon^{-s}  \| f \|_{\dot{B}^{-s}}  +\varepsilon^{\alpha}  \mathbf{G}_{\alpha,p}(f)(x)\right) \,.  \]
		Since the indicated inequality holds true for  $\varepsilon>0$, then minimizing the right hand side  of this one yields the desired result.  \\
		Hence, we complete the proof of Lemma \ref{Lem10}.
	\end{proof}
	Next, we have the following lemma.
	\begin{lemma}\label{Lem11}
		Let $0<\alpha_1<\alpha_2< 1$. Let $1\leq p_1, p_2 <\infty$, and  $r>1$ be such that 
		\begin{equation}\label{1.6}
			\frac{1}{p_1}= \frac{1}{r} \left(1-\frac{\alpha_1}{\alpha_2}\right) + \frac{1}{p_2}\frac{\alpha_1}{\alpha_2}  \,. 
		\end{equation} 
		If  $f\in L^r(\RR^n)\cap \dot{W}^{\alpha_2,p_2}(\RR^n)$, then $f\in \dot{W}^{\alpha_1,p_1}(\RR^n)$. In addition,  there exists a constant $C=C(\alpha_1,\alpha_2,p_1,p_2,n)>0$ such that
		\begin{equation}\label{1.8}
			\|f\|_{\dot{W}^{\alpha_1,p_1}}  \leq C 
			\| f \|^{1-\frac{\alpha_1}{\alpha_2}}_{ L^{r} }   
			\|f\|^{\frac{\alpha_1}{\alpha_2}}_{\dot{W}^{\alpha_2,p_2}}  \,.
		\end{equation}	
	\end{lemma}
	\begin{proof}[Proof of Lemma \ref{Lem11}] 
		For any set $\Omega$ in $\mathbb{R}^n$, let us denote  $\fint_\Omega  f(x) \, dx =  \frac{1}{|\Omega|}  \int_{\Omega}   f(x) \, dx$. 
		\\
		For any  $x, z\in\mathbb{R}^n$, we have from the triangle inequality and  change of variables that
		\begin{align*}
			\big|f(x+z)-f(x)\big|&\leq \big|f(x+z) - \fint_{B(x,|z|)}  f(y)\,dy  \big|+\big|f(x) - \fint_{B(x,|z|)}  f(y)\,dy \big|  
			\\
			&\leq \fint_{B(x,|z|)}  \big|f(x+z) -f(y) \big|  \, dy+ \fint_{B(x,|z|)}  \big|f(x)-f(y) \big|  \, dy 
			\\
			&\leq C(n)  \left( \fint_{B(0,2|z|)}  \big|f(x+z) -f(x+z+y) \big|  \, dy+ \fint_{B(0,2|z|)}  \big|f(x) - f(x+y) \big|  \, dy   \right)\,.
		\end{align*}
		With the last inequality noted, and by using a change of variables, we get
		\begin{align}\label{1.10}
			\int\int \frac{|f(x+z)-f(x)|^{p_1}}{|z|^{n+\alpha_1 p_1}} dzdx   \lesssim \int\int  \left( \fint_{B(0,2|z|)}  \big|f(x)-f(x+y) \big|  \, dy\right)^{p_1}  \frac{dzdx}{|z|^{n+\alpha_1 p_1}}\,.
		\end{align}
		Next, for every   $p\geq 1$ we show that   
		\begin{align}\label{1.13}
			\int \left( \fint_{B(0,2|z|)} \big|f(x)-f(x+y) \big|  \, dy\right)^{p_1}  \frac{dz}{|z|^{n+\alpha_1 p_1}}  \lesssim \big[\mathbf{M}(f)(x)\big]^{\frac{(\alpha_2-\alpha_1)p_1}{\alpha_2}}  \left[\mathbf{G}_{\alpha_2,p}(x) \right]^{\frac{\alpha_1 p_1}{\alpha_2}}.
		\end{align}
		Thanks to Remark \ref{Rem4},  it suffices to show that \eqref{1.13} holds for $p=1$.
		\\
		Indeed, we have
		\begin{align}\label{1.14}
			\int_{\{|z|<t\}} \left( \fint_{B(0,2|z|)}  \big|f(x) - f(x+y) \big|  \, dy\right)^{p_1}  \frac{dz}{|z|^{n+\alpha_1 p_1}} 
			&= \int_{\{|z|<t\}}  \left( \fint_{B(0,2|z|)}  \frac{\big|f(x)-f(x+y) \big|}{|z|^{\alpha_2}}  \, dy\right)^{p_1}  \frac{|z|^{\alpha_2 p_1}  dz}{|z|^{n+\alpha_1 p_1}}  \nonumber	\\
			&\lesssim  	 \left[ \mathbf{G}_{\alpha_2,1}(x)\right]^{p_1} \int_{\{|z|<t\}}  \frac{1}{|z|^{n+(\alpha_1-\alpha_2) p_1}}  dz \nonumber
			\\
			&\lesssim	t^{(\alpha_2-\alpha_1) p_1}  \left[ \mathbf{G}_{\alpha_2,1}(x)\right]^{p_1}\,.
		\end{align}
		On the other hand, it is not difficult to observe that
		\begin{align}\label{1.15}
			\int_{|z|\geq t}  \left( \fint_{B(0,2|z|)}  \big|f(x) - f(x+y) \big|  \, dy\right)^{p_1}  \frac{dz}{|z|^{n+\alpha_1 p_1}} 
			&\lesssim  \big[\mathbf{M}(f)(x)\big]^{p_1} \left( \int_{|z|\geq  t}\frac{dz}{|z|^{n+\alpha_1 p_1}}  \right) \nonumber
			\\
			&\lesssim t^{-\alpha_1 p_1}\big[\mathbf{M}(f)(x)\big]^{p_1}  \,.
		\end{align}
		From \eqref{1.14} and \eqref{1.15}, we obtain
		\begin{align*}
			\int \left( \fint_{B(0,2|z|)}  \big|f(x) - f(x+y) \big|  \, dy\right)^{p_1}  \frac{dz}{|z|^{n+\alpha_1 p_1}} \lesssim  	t^{(\alpha_2-\alpha_1) p_1}  \left[ \mathbf{G}_{\alpha_2,1}(x)\right]^{p_1}  +t^{-\alpha_1 p_1}\big[\mathbf{M}(f)(x)\big]^{p_1}  \,.
		\end{align*}
		Minimizing the right hand side of the last inequality yields \eqref{1.13}.
		\\
		Then, it follows from \eqref{1.13} that
		\[ 
		\int\int \left( \fint_{B(0,2|z|)} \big|f(x)-f(x+y) \big|  \, dy\right)^{p_1}  \frac{dzdx}{|z|^{n+\alpha_1 p_1}} dx \lesssim  \int \big[\mathbf{M}(f)(x)\big]^{\frac{(\alpha_2-\alpha_1)p_1}{\alpha_2}}  \left[\mathbf{G}_{\alpha_2,p_2}(x) \right]^{\frac{\alpha_1 p_1}{\alpha_2}}   dx\,.
		\]
		Note that $\alpha_2 p_2>\alpha_1 p_1$,  and $r=\frac{p_1p_2(\alpha_2-\alpha_1)}{\alpha_2p_2-\alpha_1p_1}$, see \eqref{1.6}. 
		Then, applying  H\"older's inequality with  $\big((\frac{\alpha_2p_2}{\alpha_1p_1})^\prime, \frac{\alpha_2p_2}{\alpha_1p_1}\big)$ to the right hand side of the last inequaltiy    yields
		\begin{align*}
			\int\int  \left( \fint_{B(0,2|z|)}  \big|f(x)-f(x+y) \big|  \, dy\right)^{p_1}  \frac{dzdx}{|z|^{n+\alpha_1 p_1}}  
			&\lesssim  \big\|\mathbf{M}(f)\big\|^{\frac{(\alpha_2-\alpha_1) p_1}{\alpha_2}}_{L^r}  \|\mathbf{G}_{\alpha_2,p_2}\|^{\frac{\alpha_1p_1}{\alpha_2}}_{L^{p_2}} \,.
		\end{align*}
		Thanks to Remark \ref{Rem4}, and by the fact that $\mathbf{M}$ maps $L^r(\mathbb{R}^n)$ into $L^r(\mathbb{R}^n)$ $r>1$, we deduce  from the last inequality that 
		\begin{align}\label{1.18}
			\int\int  \left( \fint_{B(0,2|z|)}  \big|f(x)-f(x+y) \big|  \, dy\right)^{p_1}  \frac{dzdx}{|z|^{n+\alpha_1 p_1}}  \lesssim \|f\|^{\frac{(\alpha_2-\alpha_1) p_1}{\alpha_2}}_{L^r}  \|f\|^{\frac{\alpha_1 p_1}{\alpha_2}}_{\dot{W}^{\alpha_2,p_2}}  \,.
		\end{align}
		Combining \eqref{1.10} and \eqref{1.18} yields  \eqref{1.8}.
		\\
		Hence, we obtain  Lemma \ref{Lem11}.
	\end{proof}
	Now,  we can apply Lemma \ref{Lem10} and Lemma \ref{Lem11}  alternatively to get Theorem \ref{Mainthe} for  the case $0<\alpha_2<1$. Indeed, we apply \eqref{1.1}  to  $s=\sigma$,  $\alpha=\alpha_2$, $p=p_2$. Then,  
	\begin{align}\label{1.20}
		\|f\|_{L^q} & \lesssim   
		\|f\|^{\frac{\alpha_2}{\alpha_2+\sigma}}_{\dot{B}^{-\sigma}}  \big\|\mathbf{G}^{\frac{\sigma}{\alpha_2+\sigma}}_{\alpha_2,p_2}\big\|_{L^q}  =   
		\|f\|^{\frac{\alpha_2}{\alpha_2+\sigma}}_{\dot{B}^{-\sigma}}  \big\|\mathbf{G}_{\alpha_2,p_2}\big\|^{\frac{\sigma}{\alpha_2+\sigma}}_{L^{p_2}} \leq  
		\|f\|^{\frac{\alpha_2}{\alpha_2+\sigma}}_{\dot{B}^{-\sigma}} \|f\|^{\frac{\sigma}{\alpha_2+\sigma}}_{\dot{W}^{\alpha_2,p_2}}  \,,
	\end{align}
	with $q=p_2\big(\frac{\alpha_2+\sigma}{\sigma}\big)$.
	\\
	Since $p_1=p_2\big(\frac{\alpha_2+\sigma}{\alpha_1+\sigma}\big)$, then it follows from \eqref{1.6} that $r=q>1$. 
	\\
	Next, applying Lemma \ref{Lem11} yields
	\begin{equation}\label{1.19}
		\|f\|_{\dot{W}^{\alpha_1,p_1}}  \lesssim 
		\| f \|^{1-\frac{\alpha_1}{\alpha_2}}_{ L^{r} }   
		\|f\|^{\frac{\alpha_1}{\alpha_2}}_{\dot{W}^{\alpha_2,p_2}} \lesssim   
		\|f\|^{\frac{\alpha_2-\alpha_1}{\alpha_2+\sigma}}_{\dot{B}^{-\sigma}}  \|f\|^{\frac{\alpha_1+\sigma}{\alpha_2+\sigma}}_{\dot{W}^{\alpha_2,p_2}}\,.
	\end{equation}
	Hence, we obtain  Theorem \ref{Mainthe} for the case $0\leq \alpha_1<  \alpha_2<1$, $p_i<\infty$, $i=1,2$.
	\\
	
	To end {\bf Step 1}, it remains to study the case $\alpha_2=1$, i.e:
	\begin{equation} \label{1.22a} 	\|f\|_{\dot{W}^{\alpha_1,p_1}}  \leq C 
		\|f\|^{\frac{1-\alpha_1}{1+\sigma}}_{\dot{B}^{-\sigma}}  
		\|Df\|^{\frac{\alpha_1+\sigma}{1+\sigma}}_{L^{p_2}}
		\,. 
	\end{equation}
	This can be done if we show that
	\begin{equation}\label{1.21}
		\|f\|_{\dot{W}^{\alpha_1,p_1}}  \leq C 
		\| f \|^{1-\alpha_1}_{ L^{r} }   
		\|Df\|^{\alpha_1}_{L^{p_2}}  \,,
	\end{equation}
	with $1\leq r<\infty$,  $\frac{1}{p_1}= \frac{1-\alpha_1}{r}  + \frac{\alpha_1}{p_2}$.
	\\
	Indeed, a combination of \eqref{1.21} and \eqref{1.1a} implies that 
	\begin{align*}
		\|f\|_{\dot{W}^{\alpha_1,p_1}}  \lesssim  \|f\|^{1-\alpha_1}_{L^r}  \|Df\|^{\alpha_1}_{L^{p_2}}  \lesssim  \|f\|^{\frac{1-\alpha_1}{1+\sigma}}_{\dot{B}^{-\sigma}}  
		\|D f\|^{\frac{\sigma(1-\alpha_1)}{1+\sigma}}_{L^{p_2}} \|Df\|^{\alpha_1}_{L^{p_2}}  =   \|f\|^{\frac{1-\alpha_1}{1+\sigma}}_{\dot{B}^{-\sigma}}  
		\|D f\|^{\frac{\alpha_1+\sigma}{1+\sigma}}_{L^{p_2}}\,.
	\end{align*}
	Note that  $p_1=p_2\big( \frac{1+\sigma}{\alpha_1+\sigma}\big)$, and  $r=p_2\big( \frac{1+\sigma}{\sigma}\big)$.
	\\
	Hence, we obtain  Theorem \ref{Mainthe}  when $\alpha_2=1$.  
	\\
	Now, it remains to prove \eqref{1.21}. We note that \eqref{1.21} was proved for $p_2=1$ (see,  e.g., \cite{Brezis3, CDPX}). In fact, one can modify the proofs in \cite{Brezis3, CDPX} in order to obtain \eqref{1.21} for the case $1<p_2<\infty$. However, for consistency, we give the proof of \eqref{1.21} for $1<p_2<\infty$.   
	\\
	To obtain the result, we prove a version of \eqref{1.13} in terms of  $\mathbf{M}(|Df|)(x)$ instead of $\mathbf{G}_{1,p}(x)$. Precisely, we show that
	\begin{align}\label{1.23}
		\int \left( \fint_{B(0,2|z|)} \big|f(x)-f(x+y) \big|  \, dy\right)^{p_1}  \frac{dz}{|z|^{n+\alpha_1 p_1}}  \lesssim \big[\mathbf{M}(f)(x)\big]^{(1-\alpha_1)p_1}  \left[\mathbf{M}(|Df|)(x) \right]^{\alpha_1 p_1}
	\end{align} 
	for $x\in\RR^n$.
	\\
	Indeed, it follows from the mean value theorem and a change of variables that
	\begin{align*}
		\fint_{B(0,2|z|)}  \frac{\big|f(x)-f(x+y) \big|}{|z|}  \, dy  &\lesssim	\fint_{B(0,2|z|)}  \frac{\big|f(x)-f(x+y) \big|}{|y|}  \, dy  
		\\
		&=  	\fint_{B(0,2|z|)}  \frac{\big|\int^1_0  D f(x+\tau y)  \cdot y\,  d\tau\big| }{|y|}  \, dy 
		\\
		&\leq \int^1_0\fint_{B(x,2\tau|z|)}   | D f(\zeta) |   \,  d\zeta d\tau 
		\leq  \int^1_0 \mathbf{M}(|Df|)(x)   \,  d\tau  =  \mathbf{M}(|Df|)(x) \,.
	\end{align*}
	Thus,
	\begin{align}\label{1.24}
		\int_{\{|z|<t\}} \left( \fint_{B(0,2|z|)}  \big|f(x) - f(x+y) \big|  \, dy\right)^{p_1}  \frac{dz}{|z|^{n+\alpha_1 p_1}} 
		&= \int_{\{|z|<t\}}  \left( \fint_{B(0,2|z|)}  \frac{\big|f(x)-f(x+y) \big|}{|z|}  \, dy\right)^{p_1}  \frac{|z|^{p_1}  dz}{|z|^{n+\alpha_1 p_1}}  \nonumber
		\\
		&\lesssim   \big[\mathbf{M}(|Df|)(x)  \big]^{p_1}  \int_{\{|z|<t\}} |z|^{-n+(1-\alpha_1)p_1} \, dz \nonumber
		\\
		&\lesssim  t^{(1-\alpha_1)p_1} \big[\mathbf{M}(|Df|)(x)  \big]^{p_1} \,.
	\end{align}
	From \eqref{1.24} and \eqref{1.15}, we  obtain
	\begin{align}\label{1.25}
		\int \left( \fint_{B(0,2|z|)}  \big|f(x) - f(x+y) \big|  \, dy\right)^{p_1}  \frac{dz}{|z|^{n+\alpha_1 p_1}} \lesssim  	t^{(1-\alpha_1) p_1}  \left[ \mathbf{M}(|Df|)(x)\right]^{p_1}  +t^{-\alpha_1 p_1}\big[\mathbf{M}(f)(x)\big]^{p_1}  \,.
	\end{align} 
	Hence, \eqref{1.23} follows by minimizing the right hand side of \eqref{1.25} with respect to $t$.

	If $p_2>1$,  then we apply H\"older's inequality  in \eqref{1.23} in order to get
	\begin{align*}
		\|f\|_{\dot{W}^{\alpha_1,p_1}} &\lesssim 	\int \int \left( \fint_{B(0,2|z|)} \big|f(x)-f(x+y) \big|  \, dy\right)^{p_1}  \frac{dz  dx}{|z|^{n+\alpha_1 p_1}} 
		\\
		&\lesssim \int\big[\mathbf{M}(f)(x)\big]^{(1-\alpha_1)p_1}  \left[\mathbf{M}(|Df|)(x) \right]^{\alpha_1 p_1}  dx
		\\
		&\leq \|\mathbf{M}(f)\|^{(1-\alpha_1)p_1}_{L^r}  \|\mathbf{M}(|Df|)\|^{\alpha_1 p_1}_{L^{p_2}}  
		\\
		&\lesssim \|f\|^{(1-\alpha_1)p_1}_{L^r}  \|Df\|^{\alpha_1 p_1}_{L^{p_2}}\,, 
	\end{align*} 
	where $r=p_2\big(\frac{1+\sigma}{\sigma}\big)>1$. Note that the last inequality follows from the $L^{p}$-boundedness of $\mathbf{M}$, $p>1$.
	Thus, we get \eqref{1.21}.
	\\
	This puts an end to the proof of {\bf Step 1}.
	\\
	
	{\bf  ii) Step 2.}
	Now, we can prove Theorem \ref{Mainthe} for the case $\alpha_1\geq 1$.
	At the beginning, let us denote $\alpha_i=\floor{\alpha_i}+s_i$, $i=1,2$. Then, we divide the proof into the following cases.
	\\
	
	{\bf a) The case $\floor{\alpha_2}=\floor{\alpha_1}$:} 
	By applying Theorem \ref{Mainthe} to  $D^{\floor{\alpha_1}}  f$, $\sigma_{{\rm new}}=\sigma+\floor{\alpha_1}$;  and   by  Proposition \ref{Pro1}, we obtain
	\begin{align*}
		\big\|f\big\|_{\dot{W}^{\alpha_1,p_1}}  =	\big\|D^{\floor{\alpha_1}}  f \big\|_{\dot{W}^{s_1,p_1}} &\lesssim \big\| D^{\floor{\alpha_1}}  f \big\|^{\frac{s_2-s_1}{s_2+\sigma+\floor{\alpha_1}}}_{\dot{B}^{-(\sigma+\floor{\alpha_1})}} \big\|D^{\floor{\alpha_1}}  f \big\|^{\frac{s_1+\sigma+\floor{\alpha_1}}{s_2+\sigma+\floor{\alpha_1}}}_{\dot{W}^{s_2,p_2}}  \\
		&\lesssim \big\| f \big\|^{\frac{s_2-s_1}{s_2+\sigma+\floor{\alpha_1}}}_{\dot{B}^{-\sigma}}  
		\big\| D^{\floor{\alpha_2}}  f \big\|^{\frac{s_1+\sigma+\floor{\alpha_1}}{s_2+\sigma+\floor{\alpha_1}}}_{\dot{W}^{s_2,p_2}}  =  \big\| f \big\|^{\frac{\alpha_2-\alpha_1}{\alpha_2+\sigma}}_{\dot{B}^{-\sigma}}  
		\big\| f \big\|^{\frac{\alpha_1+\sigma}{\alpha_2+\sigma}}_{\dot{W}^{\alpha_2,p_2}} \,,
	\end{align*}
	with $p_1=p_2\big(\frac{s_2+\sigma_{{\rm new}}}{s_1+\sigma_{{\rm new}}}\big)=p_2\big(\frac{\alpha_2+\sigma}{\alpha_1+\sigma}\big)$.
	\\
	Hence, we get the conclusion for this case.
	\\
	
	{\bf b) The case $\floor{\alpha_2}>\floor{\alpha_1}$:}  If $s_2>0$, then we can apply Theorem \ref{Mainthe} to $D^{\floor{\alpha_2}}  f$,  $\sigma_{{\rm new}}=\sigma+\floor{\alpha_2}$.   Therefore,
	\begin{align}\label{1.30}
		\big\| D^{\floor{\alpha_2}}  f \big\|_{L^q} \lesssim \big\|D^{\floor{\alpha_2}}  f \big\|^{\frac{s_2}{s_2+\sigma+\floor{\alpha_2}}}_{\dot{B}^{-(\sigma+\floor{\alpha_2})}}  \big\| D^{\floor{\alpha_2}}  f \big\|^{\frac{\sigma+\floor{\alpha_2}}{s_2+\sigma+\floor{\alpha_2}}}_{\dot{W}^{s_2,p_2}} \lesssim  \big\| f \big\|^{\frac{s_2}{\alpha_2+\sigma}}_{\dot{B}^{-\sigma}}  \big\|f \big\|^{\frac{\floor{\alpha_2}+\sigma}{\alpha_2+\sigma}}_{\dot{W}^{\alpha_2,p_2}} \,,
	\end{align}
	with   $q=p_2\big(\frac{\alpha_2+\sigma}{\floor{\alpha_2}+\sigma}\big)$. Again, the last inequality follows from the lifting property  in Proposition \eqref{Pro1}.  
	\\
	Next, applying Theorem \ref{Mainthe} to $D^{\floor{\alpha_1}} f$, $\sigma_{{\rm new}}=\sigma+\floor{\alpha_1}$  yields
	\begin{align}\label{1.31}
		\big\|f\big\|_{\dot{W}^{\alpha_1,p_1}}  = 	\big\|D^{\floor{\alpha_1}}  f \big\|_{\dot{W}^{s_1,p_1}} &\lesssim \big\|D^{\floor{\alpha_1}}  f \big\|^{\frac{1-s_1}{1+\sigma+\floor{\alpha_1}}}_{\dot{B}^{-(\sigma+\floor{\alpha_1})}}  \big\| D^{\floor{\alpha_1}+1}  f \big\|^{\frac{s_1+\sigma+\floor{\alpha_1}}{1+\sigma+\floor{\alpha_1}}}_{L^{q_1}}  \nonumber
		\\
		&\lesssim \big\| f \big\|^{\frac{1-s_1}{1+\sigma+\floor{\alpha_1}}}_{\dot{B}^{-\sigma}}  \big\| D^{\floor{\alpha_1}+1}  f \big\|^{\frac{s_1+\sigma+\floor{\alpha_1}}{1+\sigma+\floor{\alpha_1}}}_{L^{q_1}} \,, 
	\end{align}
	with $q_1= p_1\big(\frac{s_1+\sigma+\floor{\alpha_1}}{1+\sigma+\floor{\alpha_1}}\big)$.
	\\
	
	If $\floor{\alpha_2}=\floor{\alpha_1}+1$, then observe that $q=q_1$. Thus, we deduce from \eqref{1.30} and \eqref{1.31} that
	\begin{align*}
		\big\|f\big\|_{\dot{W}^{\alpha_1,p_1}} \lesssim  \big\| f \big\|^{\frac{1-s_1}{1+\sigma+\floor{\alpha_1}}}_{\dot{B}^{-\sigma}}  
		\left( \big\| f \big\|^{\frac{s_2}{\alpha_2+\sigma}}_{\dot{B}^{-\sigma}}  \big\| f \big\|^{\frac{\floor{\alpha_2}+\sigma}{\alpha_2+\sigma}}_{\dot{W}^{\alpha_2,p_2}}  \right)^{\frac{\alpha_1+\sigma}{\floor{\alpha_2}+\sigma}}=  \big\| f \big\|^{\frac{\alpha_2-\alpha_1}{\alpha_2+\sigma}}_{\dot{B}^{-\sigma}}  \big\| f \big\|^{\frac{\alpha_1+\sigma}{\alpha_2+\sigma}}_{\dot{W}^{\alpha_2,p_2}} \,.
	\end{align*}
	This yields \eqref{-3}.
	\\
	Note that $\frac{1-s_1}{1+\sigma+\floor{\alpha_1}}+\frac{s_2(\alpha_1+\sigma)}{(\alpha_2+\sigma)(\floor{\alpha_2}+\sigma)}=\frac{\alpha_2-\alpha_1}{\alpha_2+\sigma}$  since $\floor{\alpha_2}=\floor{\alpha_1}+1$.
	\\
	
	If $\floor{\alpha_2}>\floor{\alpha_1}+1$, then we apply \cite[Theorem 1.2]{DaoLamLu1} to $k=\floor{\alpha_1}+1$, $m=\floor{\alpha_2}$. Thus,  
	\begin{align}\label{1.32}
		\big\| D^{\floor{\alpha_1}+1}  f \big\|_{L^{q_1}}\lesssim \big\|f\big\|^{\frac{\floor{\alpha_2}-\floor{\alpha_1}-1}{\floor{\alpha_2}+\sigma}}_{\dot{B}^{-\sigma}}  	\big\| D^{\floor{\alpha_2}}  f \big\|^{\frac{\floor{\alpha_1}+1+\sigma}{\floor{\alpha_2}+\sigma}}_{L^{q_2}}\,,
	\end{align}
	with $q_2=q_1 \big(\frac{\floor{\alpha_1}+1+\sigma}{\floor{\alpha_2}+\sigma}\big)$.
	\\
	Combining \eqref{1.31} and \eqref{1.32} yields
	\begin{align}\label{1.33}
		\big\|f \big\|_{\dot{W}^{\alpha_1,p_1}} 
		&\lesssim \big\| f \big\|^{\frac{1-s_1}{1+\sigma+\floor{\alpha_1}}}_{\dot{B}^{-\sigma}}  \left( \big\|f\big\|^{\frac{\floor{\alpha_2}-\floor{\alpha_1}-1}{\floor{\alpha_2}+\sigma}}_{\dot{B}^{-\sigma}} 
		\big\|D^{\floor{\alpha_2}}  f \big\|^{\frac{\floor{\alpha_1}+1+\sigma}{\floor{\alpha_2}+\sigma}}_{L^{q_2}}\right)^{\frac{\alpha_1+\sigma}{1+\floor{\alpha_1}+\sigma}}  \nonumber
		\\
		&=	 \big\| f \big\|^{\frac{1-s_1}{1+\sigma+\floor{\alpha_1}}  +  \big(\frac{\floor{\alpha_2}-\floor{\alpha_1}-1}{\floor{\alpha_2}+\sigma}\big)\big( \frac{\alpha_1+\sigma}{\floor{\alpha_1}+1+\sigma}\big)}_{\dot{B}^{-\sigma}}  	
		\big\| D^{\floor{\alpha_2}}  f \big\|^{\frac{\alpha_1+\sigma}{\floor{\alpha_2}+\sigma}}_{L^{q_2}}  \,.
	\end{align}
	Observe that $q=q_2=p_2 \big(\frac{\alpha_2+\sigma}{\floor{\alpha_2}+\sigma}\big) $. Thus, it follows from \eqref{1.33} and \eqref{1.30} that
	\begin{align*}
		\big\|  f \big\|_{\dot{W}^{\alpha_1,p_1}} 
		&\lesssim  \big\| f \big\|^{\frac{1-s_1}{1+\sigma+\floor{\alpha_1}}  +  \big(\frac{\floor{\alpha_2}-\floor{\alpha_1}-1}{\floor{\alpha_2}+\sigma}\big)\big( \frac{\alpha_1+\sigma}{\floor{\alpha_1}+1+\sigma}\big)}_{\dot{B}^{-\sigma}}  \left( \big\| f \big\|^{\frac{s_2}{\alpha_2+\sigma}}_{\dot{B}^{-\sigma}} \big\| f \big\|^{\frac{\floor{\alpha_2}+\sigma}{\alpha_2+\sigma}}_{\dot{W}^{\alpha_2,p_2}}\right)^{\frac{\alpha_1+\sigma}{\floor{\alpha_2}+\sigma}}
		\\
		&=\big\| f \big\|^{\frac{\alpha_2-\alpha_1}{\alpha_2+\sigma}}_{\dot{B}^{-\sigma}} \big\| f \big\|^{\frac{\alpha_1+\sigma}{\alpha_2+\sigma}}_{\dot{W}^{\alpha_2,p_2}}   \,.
	\end{align*} 
	A straightforward computation shows that 
	\[\frac{1-s_1}{1+\sigma+\floor{\alpha_1}}  +  \big(\frac{\floor{\alpha_2}-\floor{\alpha_1}-1}{\floor{\alpha_2}+\sigma}\big)\big( \frac{\alpha_1+\sigma}{\floor{\alpha_1}+1+\sigma}\big)+ \frac{s_2(\alpha_1+\sigma)}{(\alpha_2+\sigma)(\floor{\alpha_2}+\sigma)} =  \frac{\alpha_2-\alpha_1}{\alpha_2+\sigma} \,.\]
	This puts an end to the proof of Theorem \ref{Mainthe} for  $s_2>0$.
	\\
	The proof of the case $s_2=0$ can be done similarly as above. Then, we leave the details to the reader.
	\\
	Hence, we complete the proof of Theorem \ref{Mainthe}. 

	\subsection{Proof of Theorem \ref{Mainthe1}}At the beginning, let us recall the notation  $\alpha_i=\floor{\alpha_i}+s_i$, $i=1, 2$. Then,	we divide the proof into the two following cases.
	\\
	
	{\bf i) The case $p_1=p_2=\infty$.} If  $0<\alpha_1<\alpha_2<1$, then  \eqref{4.1} becomes
	\begin{equation}\label{2.0}
		\|f\|_{\dot{C}^{\alpha_1}}  \lesssim
		\|f\|^{1-\frac{\alpha_1}{\alpha_2}}_{\dot{B}^0} \|f\|^{\frac{\alpha_1}{\alpha_2}}_{\dot{C}^{\alpha_2}}  \,.
	\end{equation}
	Inequality \eqref{2.0} can be obtained easily from  the proof of \eqref{1.0} for $\sigma=0$. Then, we leave the details to the reader.
	\\
	
	If $0<\alpha_1<1\leq \alpha_2$, and $\alpha_2$ is integer,  then   \eqref{4.1} reads as:
	\begin{equation}\label{2.0a}
		\|f\|_{\dot{C}^{\alpha_1}}  \lesssim
		\|f\|^{1-\frac{\alpha_1}{\alpha_2}}_{\dot{B}^0} \|\nabla^{\alpha_2} f\|^{\frac{\alpha_1}{\alpha_2}}_{L^\infty}  \,.
	\end{equation}
	To obtain \eqref{2.0a}, we utilize  the vanishing moments of  $\varphi_\varepsilon$  in Theorem \ref{ThePeetre}. In fact, let us fix $k>\alpha_2$. 
	Then, it follows from the Taylor series  that
	\begin{align}\label{2.0c}
		|\varphi_\varepsilon * f(x)| &=\left|\int\big( f(x-y)-f(x)\big) \varphi_\varepsilon(y)\, dy   \right|  \nonumber 
		\\
		&=\left| \int \left(\sum_{|\gamma|< \alpha_2} \frac{D^\gamma f(x)}{|\gamma|!}  (-y)^\gamma + \sum_{|\gamma|= \alpha_2} \frac{D^\gamma f(\zeta)}{|\gamma|!}  (-y)^\gamma  \right) \varphi_\varepsilon(y)\,dy \right|
		\nonumber	\\
		&=\left|  \int \sum_{|\gamma|= \alpha_2} \frac{D^{\gamma} f(\zeta)}{|\alpha_2|!}  (-y)^\gamma   \varphi_\varepsilon(y)\,dy \right|
	\end{align}
	for some $\zeta$ in the line-$xy$. Note that $$ \int  \frac{D^\gamma f(x)}{|\gamma|!}  (-y)^\gamma \varphi_\varepsilon(y)\,dy=0$$ for every multi-index $|\gamma|<k$. 
	\\
	Hence, we get from \eqref{2.0c} that
	\[  |\varphi_\varepsilon * f(x)|	 \lesssim   \|\nabla^{\alpha_2} f\|_{L^\infty}   \int_{B(0,\varepsilon)}   |y|^{\alpha_2}  |\varphi_\varepsilon(y)|\,dy\lesssim   \varepsilon^{\alpha_2}  \|\nabla^{\alpha_2} f\|_{L^\infty}\,.\]
	Inserting  the last inequality into \eqref{2.-1} yields
	\[ 	\varepsilon^{-\alpha_1}\|\varphi_\varepsilon  * f\|_{L^\infty}  \lesssim   \delta^{\alpha_2-\alpha_1} \|\nabla^{\alpha_2} f\|_{L^\infty}
	+\delta^{-\alpha_1} \|f\|_{\dot{B}^{0}}  \,.   
	\] 
	By minimizing the  right hand side of the indicated inequality, we get
	\[ 	
	\varepsilon^{-\alpha_1}\|\varphi_\varepsilon  * f\|_{L^\infty}  \lesssim 	\|f\|^{1-\frac{\alpha_1}{\alpha_2}}_{\dot{B}^0} \|\nabla^{\alpha_2} f\|^{\frac{\alpha_1}{\alpha_2}}_{L^\infty} \,.
	\]
	This implies \eqref{2.0a}.
	\\
	
	If $0<\alpha_1<1\leq \alpha_2$, and $\alpha_2$ is not integer,  then   \eqref{4.1} reads as:
	\begin{equation}\label{2.0b}
		\|f\|_{\dot{C}^{\alpha_1}}  \lesssim
		\|f\|^{1-\frac{\alpha_1}{\alpha_2}}_{\dot{B}^0} \left\|\nabla^{\floor{\alpha_2}} f\right\|^{\frac{\alpha_1}{\alpha_2}}_{\dot{C}^{s_2}}  \,.
	\end{equation}
	To obtain \eqref{2.0b}, we apply \eqref{2.0c} to $\floor{\alpha_2}$. Thus, 
	\begin{align*}
		|\varphi_\varepsilon * f(x)| 	&=\left|  \int \sum_{|\gamma|=\floor{\alpha_2}} \frac{D^{\gamma} f(\zeta)}{\floor{\alpha_2}!}  (-y)^\gamma   \varphi_\varepsilon(y)\,dy \right|
		\\
		&=\left|  \int \sum_{|\gamma|= \floor{\alpha_2}} \frac{D^{\gamma} f(\zeta) -  D^{\gamma} f(x)}{\floor{\alpha_2}!}  (-y)^\gamma   \varphi_\varepsilon(y)\,dy \right|
		\\
		&\lesssim \|D^{\floor{\alpha_2}} f\|_{\dot{C}^{s_2}} \int   |x-\zeta|^{s_2} |y|^{\floor{\alpha_2}}  |\varphi_\varepsilon(y)|\,dy 
		\\
		&\leq\|D^{\floor{\alpha_2}} f\|_{\dot{C}^{s_2}} \int_{B(0,\varepsilon)}   |y|^{s_2} |y|^{\floor{\alpha_2}}  |\varphi_\varepsilon(y)|\,dy \lesssim \varepsilon^{\alpha_2} 
		\|D^{\floor{\alpha_2}} f\|_{\dot{C}^{s_2}}  \,.
	\end{align*}
	Thus, 
	\[ 	\varepsilon^{-\alpha_1}\|\varphi_\varepsilon  * f\|_{L^\infty}  \lesssim   \delta^{\alpha_2-\alpha_1} \|D^{\floor{\alpha_2}} f\|_{\dot{C}^{s_2}}
	+\delta^{-\alpha_1} \|f\|_{\dot{B}^{0}}  \,.   \] 
	By the analogue as in the proof of \eqref{2.0a}, we also obtain \eqref{2.0b}.  
	\\
	In conclusion, Theorem \ref{Mainthe1} was proved for the case $0<\alpha_1<1$.
	\\
	
	Now, if $\alpha_1\geq 1$,  then \eqref{4.1} becomes
	\begin{equation}\label{2.0d}
		\| D^{\floor{\alpha_1}} f \|_{\dot{C}^{s_1}}  \lesssim
		\|f\|^{1-\frac{\alpha_1}{\alpha_2}}_{\dot{B}^0} \|D^{\floor{\alpha_2}} f\|^{\frac{\alpha_1}{\alpha_2}}_{\dot{C}^{s_2}}  \,.
	\end{equation}
	Again, we  note that $\|.\|_{\dot{C}^{s_i}}$ is replaced by $\|.\|_{L^\infty}$ whenever $s_i=0$, $i=1,2$. 
	\\
	To obtain  \eqref{2.0d}, we apply Theorem \ref{Mainthe} to $f_{{\rm new}}=D^{\floor{\alpha_1}} f$, and $\sigma =\floor{\alpha_1}$.  
	\\
	Hence, it follows from Proposition \ref{Pro1} that
	\begin{align*}
		\|f\|_{\dot{C}^{\alpha_1}} =\| D^{\floor{\alpha_1}} f \|_{\dot{C}^{s_1}}  
		&\lesssim
		\big\|D^{\floor{\alpha_1}} f\big\|^{\frac{\alpha_2-\floor{\alpha_1}-s_1}{\alpha_2-\floor{\alpha_1}+\sigma}}_{\dot{B}^{-\floor{\alpha_1}}} \big\|D^{\floor{\alpha_1}} f\big\|^{\frac{s_1+\sigma}{\alpha_2-\floor{\alpha_1}+\sigma}}_{\dot{C}^{\alpha_2-\floor{\alpha_1}}}
		\\
		&\lesssim
		\|f\|^{\frac{\alpha_2-\alpha_1}{\alpha_2}}_{\dot{B}^{0}} 
		\big\|D^{\floor{\alpha_2}}f\big\|^{\frac{\alpha_1}{\alpha_2}}_{\dot{C}^{s_2}} =	\|f\|^{\frac{\alpha_2-\alpha_1}{\alpha_2}}_{\dot{B}^{0}}
		\|f\|^{\frac{\alpha_1}{\alpha_2}}_{\dot{C}^{\alpha_2}}  \,.
	\end{align*}
	This puts an end to the proof of Theorem \ref{Mainthe1} for the case  $p_1=p_2=\infty$.
	\\
	
	{\bf ii) The case $p_i<\infty, i=1,2$.}
	We first consider the case $0<\alpha_1<1$. 
	\\
	
	{\bf a)}	If $\alpha_2\in(\alpha_1, 1)$,  then we utilize the following result  $\|\cdot\|_{\dot{W}^{s,p}}  \approx  \|\cdot\|_{\dot{B}^{s}_{p,p}}$ for $s\in(0,1)$, $p\geq 1$, see Proposition \ref{Pro-cha} in the Appendix section. Therefore, \eqref{4.1} is equivalent to the following inequality
	\begin{equation}\label{2.1}
		\|f\|_{\dot{B}^{\alpha_1}_{p_1,p_1}}  \lesssim \|f\|^{1-\frac{\alpha_1}{\alpha_2}}_{\dot{B}^{0}}\|f\|^\frac{\alpha_1}{\alpha_2}_{\dot{B}^{\alpha_2}_{p_2,p_2}}  \,.
	\end{equation}
	Note that $\alpha_1 p_1 = \alpha_2 p_2$. Hence,
	\begin{align}\label{2.2}
		2^{jn \alpha_1 p_1} \|f*\phi_j\|_{L^{p_1}}^{p_1} \leq 	2^{jn \alpha_2 p_2} \|f*\phi_j\|_{L^{p_2}}^{p_2} \|f*\phi_j\|_{L^{\infty}}^{p_1-p_2}  \leq  	2^{jn \alpha_2 p_2} \|f*\phi_j\|_{L^{p_2}}^{p_2}  \|f\|_{\dot{B}^0}^{p_1-p_2}  \,.
	\end{align} 
	This implies that 
	\[ 	\|f\|^{p_1}_{\dot{B}^{\alpha_1}_{p_1,p_1}}  \leq \|f\|^{p_1-p_2}_{\dot{B}^{0}}  \|f\|^{p_2}_{\dot{B}^{\alpha_2}_{p_2,p_2}} \]
	which is \eqref{2.1}.
	\\
	
	{\bf  b)}  If  $\alpha_2= 1$, then  we  show that
	\begin{equation}\label{2.3}
		\|f\|_{\dot{W}^{\alpha_1,p_1}}  \lesssim\|f\|^{1-\alpha_1}_{\dot{B}^{0}}  \|Df\|^{\alpha_1}_{L^{p_2}}\,. 
	\end{equation}
	To obtain \eqref{2.3}, we prove the homogeneous version of \eqref{-16}.
	\begin{lemma}\label{Lem-Hom-Sobolev}
		Let  $0<\alpha_0<\alpha_1 <\alpha_2\leq 1$, and $p_0\geq 1$ be such that  $\alpha_0 -\frac{1}{p_0}<\alpha_2-\frac{1}{p_2}$, and 
		\[  
		\frac{1}{p_1} = \frac{\theta}{p_0}   +  \frac{1-\theta}{p_2}  ,\quad  \theta= \frac{\alpha_2-\alpha_1}{\alpha_2-\alpha_0} 
		\,. \]
		Then,  we have
		\begin{equation}\label{2.4}
			\|f\|_{\dot{W}^{\alpha_1,p_1}}  \lesssim 	\|f\|^{\frac{\alpha_2-\alpha_1}{\alpha_2-\alpha_0} }_{\dot{W}^{\alpha_0,p_0}}  	\|f\|^{\frac{\alpha_1-\alpha_0}{\alpha_2-\alpha_0}}_{\dot{W}^{\alpha_2,p_2}},\quad \forall f\in   \dot{W}^{\alpha_0,p_0} (\RR^n) \cap \dot{W}^{\alpha_2,p_2}(\RR^n)  \,.
		\end{equation}
	\end{lemma}
	\begin{proof}[Proof of Lemma \ref{Lem-Hom-Sobolev}]  The proof is quite similar to the one in Lemma \ref{Lem10}.
		Indeed,	the proof follows by way of the following result.
		\\
		If $f\in \dot{W}^{\alpha_0,p_0} (\RR^n) \cap \dot{W}^{\alpha_2,p_2}(\RR^n)$, then there hold true
		\begin{align}\label{2.5a}
			\int \left( \fint_{B(0,2|z|)} \big|f(x)-f(x+y) \big|  \, dy\right)^{p_1}  \frac{dz}{|z|^{n+\alpha_1 p_1}} 
			\lesssim\big[\mathbf{G}_{\alpha_0,p_0}(f)(x)\big]^{(\frac{\alpha_2-\alpha_1}{\alpha_2-\alpha_0})p_1}   \big[\mathbf{G}_{\alpha_2,p_2}(f)(x) \big]^{(\frac{\alpha_1-\alpha_0}{\alpha_2-\alpha_0})p_1}
		\end{align} 
		if provided that $\alpha_2<1$, and 
		\begin{align}\label{2.5}
			\int \left( \fint_{B(0,2|z|)} \big|f(x)-f(x+y) \big|  \, dy\right)^{p_1}  \frac{dz}{|z|^{n+\alpha_1 p_1}} 
			\lesssim\big[\mathbf{G}_{\alpha_0,p_0}(f)(x)\big]^{(\frac{1-\alpha_1}{1-\alpha_0})p_1}   \left[\mathbf{M}(|Df|)(x) \right]^{(\frac{\alpha_1-\alpha_0}{1-\alpha_0})p_1}
		\end{align} 
		if $\alpha_2=1$.
		\\
		The proof of \eqref{2.5a} (resp. \eqref{2.5}) can be done  similarly as  the one of \eqref{1.13} (resp. \eqref{1.23}).  Therefore,  we only need to replace   $\mathbf{M}(f)(x)$ by    $\mathbf{G}_{\alpha_0,p_0}(f)(x)$   in \eqref{1.13} (resp. \eqref{1.23}). 
		\\
		In fact, we have from H\"{o}lder's inequality
		\begin{align}\label{2.6}
			\int_{\{|z|\geq t\}}  \left( \fint_{B(0,2|z|)}  \big|f(x) - f(x+y) \big|  \, dy\right)^{p_1}  \frac{dz}{|z|^{n+\alpha_1 p_1}} 
			&\leq   \int_{\{|z|\geq t\}}  \left( \fint_{B(0,2|z|)}  \big|f(x) - f(x+y) \big|^{p_0}  \, dy\right)^{\frac{p_1}{p_0}}  \frac{dz}{|z|^{n+\alpha_1 p_1}}    \nonumber
			\\
			&=  \int_{\{|z|\geq t\}}  \left( \fint_{B(0,2|z|)}  \frac{\big|f(x) - f(x+y) \big|^{p_0}}{|z|^{\alpha_0 p_0}}  \, dy\right)^{\frac{p_1}{p_0}}  \frac{|z|^{\alpha_0 p_1}dz}{|z|^{n+\alpha_1 p_1}}     \nonumber
			\\
			&\lesssim   \big[\mathbf{G}_{\alpha_0,p_0}(f)(x)\big]^{p_1}  \int_{\{|z|\geq t\}} |z|^{-n-(\alpha_1-\alpha_0)p_1}  \, dz\nonumber
			\\
			&\lesssim    t^{-(\alpha_1-\alpha_0)p_1} \big[\mathbf{G}_{\alpha_0,p_0}(f)(x)\big]^{p_1} \,.
		\end{align}	
		
		If $\alpha_2<1$, then it follows from \eqref{2.5} and \eqref{1.14} that
		\begin{align*}
			\int_{\RR^n}  \left( \fint_{B(0,2|z|)}  \big|f(x) - f(x+y) \big|  \, dy\right)^{p_1}  \frac{dz}{|z|^{n+\alpha_1 p_1}} \lesssim  t^{-(\alpha_1-\alpha_0)p_1} \big[\mathbf{G}_{\alpha_0,p_0}(f)(x)\big]^{p_1}+  t^{(\alpha_2-\alpha_1)p_1} \big[\mathbf{G}_{\alpha_2,p_2}(f)(x)\big]^{p_1}  \,.
		\end{align*}
		Thus,   \eqref{2.5a} follows by minimizing the right hand side of the indicated inequality.
		\\
		Next, applying H\"older's inequality  in \eqref{2.5a}  with $\big(\frac{p_0(\alpha_2-\alpha_0)}{p_1(\alpha_2-\alpha_1)},\frac{p_2(\alpha_2-\alpha_0)}{p_1(\alpha_2-\alpha_1)}\big)$  yields
		\begin{align*}
			\|f\|^{p_1}_{\dot{W}^{\alpha_1,p_1}}  &\lesssim 	\int_{\RR^n}  \int_{\RR^n} \left( \fint_{B(0,2|z|)}  \big|f(x) - f(x+y) \big|  \, dy\right)^{p_1}  \frac{dz}{|z|^{n+\alpha_1 p_1}}
			\\ 
			&\lesssim   \int  \big[\mathbf{G}_{\alpha_0,p_0}(f)(x)\big]^{(\frac{\alpha_2-\alpha_1}{\alpha_2-\alpha_0})p_1}   \big[\mathbf{G}_{\alpha_2,p_2}(f)(x)\big]^{(\frac{\alpha_1-\alpha_0}{\alpha_2-\alpha_0})p_1}  dx
			\\
			&\leq    \big\| \mathbf{G}_{\alpha_0,p_0} \big\|^{(\frac{\alpha_2-\alpha_1}{\alpha_2-\alpha_0})p_1}_{L^{p_0}}  \big\| \mathbf{G}_{\alpha_2,p_2} \big\|^{(\frac{\alpha_1-\alpha_0}{\alpha_2-\alpha_0})p_1}_{L^{p_2}}
			\\
			&\leq \|f\|^{(\frac{\alpha_2-\alpha_1}{\alpha_2-\alpha_0})p_1}_{\dot{W}^{\alpha_0,p_0}} \|f\|^{(\frac{\alpha_1-\alpha_0}{\alpha_2-\alpha_0})p_1}_{\dot{W}^{\alpha_2,p_2}}   \,.
		\end{align*}
		Note that the last inequality is obtained by Remark \ref{Rem4}. Hence, we get \eqref{2.4} for $\alpha_2<1$.
		\\
		
		If $\alpha_2=1$, then
		it follows from \eqref{2.6} and \eqref{1.24} that
		\begin{align*}
			\int_{\RR^n}  \left( \fint_{B(0,2|z|)}  \big|f(x) - f(x+y) \big|  \, dy\right)^{p_1}  \frac{dz}{|z|^{n+\alpha_1 p_1}} \lesssim  t^{-(\alpha_1-\alpha_0)p_1} \big[\mathbf{G}_{\alpha_0,p_0}(f)(x)\big]^{p_1}+  t^{(1-\alpha_1)p_1} \left[\mathbf{M}(|Df|)(x) \right]^{p_1}
		\end{align*}
		which implies \eqref{2.5}.
		\\
		By applying H\"older's inequality  with $\big(\frac{p_0(1-\alpha_0)}{p_1(1-\alpha_1)},\frac{p_2(1-\alpha_0)}{p_1(1-\alpha_1)}\big)$, we obtain
		\begin{align*}
			\|f\|^{p_1}_{\dot{W}^{\alpha_1,p_1}}  &\lesssim 	\int_{\RR^n} \int_{\RR^n} \left( \fint_{B(0,2|z|)}  \big|f(x) - f(x+y) \big|  \, dy\right)^{p_1}  \frac{dz}{|z|^{n+\alpha_1 p_1}} \\
			&\lesssim  \int_{\RR^n} \big[\mathbf{G}_{\alpha_0,p_0}(f)(x)\big]^{(\frac{1-\alpha_1}{1-\alpha_0})p_1}   \left[\mathbf{M}(|Df|)(x) \right]^{(\frac{\alpha_1-\alpha_0}{1-\alpha_0})p_1}  dx
			\\
			&\leq  \big\|\mathbf{G}_{\alpha_0,p_0}(f)\big\|^{(\frac{1-\alpha_1}{1-\alpha_0})p_1}_{L^{p_0}}   \big\|\mathbf{M}(|Df|) \big\|^{(\frac{\alpha_1-\alpha_0}{1-\alpha_0})p_1}_{L^{p_2}}  
			\\
			&\lesssim  \|f\|^{(\frac{1-\alpha_1}{1-\alpha_0})p_1}_{\dot{W}^{\alpha_0, p_0}} \|Df\|^{(\frac{\alpha_1-\alpha_0}{1-\alpha_0})p_1}_{L^{p_2}}  	\,.
		\end{align*} 
		This yields \eqref{2.4} for $\alpha_2=1$.
		\\
		Hence, we complete the proof of Lemma \ref{Lem-Hom-Sobolev}.  
	\end{proof}
	Now, we apply Lemma \ref{Lem-Hom-Sobolev} when $\alpha_2=1$ in order to obtain
	\begin{align*}
		\|f\|_{\dot{W}^{\alpha_1, p_1}}\lesssim  \|f\|^{\frac{1-\alpha_1}{1-\alpha_0}}_{\dot{W}^{\alpha_0, p_0}} \|Df\|^{\frac{\alpha_1-\alpha_0}{1-\alpha_0}}_{L^{p_2}}  	\,,
	\end{align*}  
	where $\alpha_0, p_0$ are chosen as in Lemma \ref{Lem-Hom-Sobolev}. 
	\\
	After that, we have  from \eqref{2.1} that
	\[    \|f\|_{\dot{W}^{\alpha_0, p_0}}  \lesssim \|f\|^{1-\frac{\alpha_0}{\alpha_1}}_{\dot{B}^0}	\|f\|^{\frac{\alpha_0}{\alpha_1}}_{\dot{W}^{\alpha_1, p_1}}   \,.\] 
	Combining the last two inequalities yields the desired result. 
	\\
	
	{\bf  The case $\alpha_2>1$}.
	\\
	If $\alpha_2$ is not integer, then we apply Theorem \ref{Mainthe}  to $\sigma=\floor{\alpha_2}$ to get
	\begin{align}\label{2.7}
		\| D^{\floor{\alpha_2}} f \|_{L^q}  \lesssim 
		\| D^{\floor{\alpha_2}} f \|^{\frac{s_2}{\alpha_2}}_{\dot{B}^{-\floor{\alpha_2}}}  \|D^{\floor{\alpha_2}}  f\|^{\frac{\floor{\alpha_2}}{\alpha_2}}_{\dot{W}^{s_2,p_2}}  \lesssim 
		\|f\|^{\frac{s_2}{\alpha_2}}_{\dot{B}^0}  
		\|f\|^{\frac{\floor{\alpha_2}}{\alpha_2}}_{\dot{W}^{\alpha_2,p_2}}  \,,
	\end{align}
	with $q=p_2\frac{\alpha_2}{\floor{\alpha_2}}$.
	Recall that $\alpha_2=\floor{\alpha_2}+s_2$.
	
	If $\floor{\alpha_2}=1$, then it follows from \eqref{2.3} and the last inequality that
	\begin{align*}
		\|f\|_{\dot{W}^{\alpha_1,p_1}}\lesssim \|f\|^{1-\alpha_1}_{\dot{B}^0} \|Df\|^{\alpha_1}_{L^q}  \lesssim  \|f\|^{1-\alpha_1}_{\dot{B}^0}  \left( \|f\|^{\frac{s_2}{\alpha_2}}_{\dot{B}^0}  
		\|f\|^{\frac{1}{\alpha_2}}_{\dot{W}^{\alpha_2,p_2}}  \right)^{\alpha_1}=\|f\|^{\frac{\alpha_2-\alpha_1}{\alpha_2}}_{\dot{B}^0}\|f\|^{\frac{\alpha_1}{\alpha_2}}_{\dot{W}^{\alpha_2,p_2}}  \,,
	\end{align*}
	with $q=\alpha_1  p_1=\alpha_2  p_2$ since $\floor{\alpha_2}=1$.
	\\
	This yields \eqref{4.1} when $\floor{\alpha_2}=1$.
	
	If $\floor{\alpha_2}>1$, then we can apply Theorem \ref{TheC} in order to get
	\begin{align*}
		\|Df\|_{L^{q_1}}  \lesssim    \|f\|^{\frac{\floor{\alpha_2}-1}{\floor{\alpha_2}}}_{\dot{B}^0} \big\|D^{\floor{\alpha_2}}f \big\|^{\frac{1}{\floor{\alpha_2}}}_{L^{q_2}}  \,,
	\end{align*}
	with $q_1=\alpha_1 p_1$, and $q_2 = \frac{q_1}{\floor{\alpha_2}}  =  \frac{\alpha_2 p_2}{\floor{\alpha_2}}$.
	\\
	A combination of the last inequality, \eqref{2.7}, and \eqref{2.3} implies that
	\begin{align*}
		\|f\|_{\dot{W}^{\alpha_1,p_1}}  &\lesssim \|f\|^{1-\alpha_1}_{\dot{B}^0} \|Df\|^{\alpha_1}_{L^{q_1}}  \lesssim  \|f\|^{1-\alpha_1}_{\dot{B}^0}  \left(\|f\|^{\frac{\floor{\alpha_2}-1}{\floor{\alpha_2}}}_{\dot{B}^0} \big\|D^{\floor{\alpha_2}}f \big\|^{\frac{1}{\floor{\alpha_2}}}_{L^{q_2}} \right)^{\alpha_1}
		\\
		&  \lesssim\|f\|^{1-\frac{\alpha_1}{\floor{\alpha_2}}}_{\dot{B}^0}  \left(\|f\|^{\frac{s_2}{\alpha_2}}_{\dot{B}^0}  
		\|f\|^{\frac{\floor{\alpha_2}}{\alpha_2}}_{\dot{W}^{\alpha_2,p_2}}  \right)^{\frac{\alpha_1}{\floor{\alpha_2}}}
		=  \|f\|^{1-\frac{\alpha_1}{\alpha_2}}_{\dot{B}^0}  \|f\|^{\frac{\alpha_1}{\alpha_2}}_{\dot{W}^{\alpha_2,p_2}}   \,.
	\end{align*}
	Hence, we obtain\eqref{4.1} when $\floor{\alpha_2}>1$.
	\\
	The case where $\alpha_2>1$ is integer can be done similarly as the above. Then, we leave the details to the reader.
	

	\section{Appendix}
	\begin{proposition}\label{Pro-cha}  The following statement holds true 
		\begin{equation}\label{5.1}
			\|f\|_{\dot{W}^{\alpha,p}} \approx  \|f\|_{\dot{B}^{\alpha}_{p,p}} ,\quad \forall  f\in \mathcal{S}(\RR^n) \,.	
		\end{equation}
	\end{proposition}
	\begin{proof}[Proof of Proposition \ref{Pro-cha}] To obtain the result,  we follow the proof by Grevholm, \cite{Grevholm}.
		\\  
		First of all, for any $s\in(0,1)$, $1\leq p<\infty$,	it is known that   (see,  e.g.,  \cite{Leoni, Triebel})
		\[
		\|f\|_{\dot{W}^{s,p}} \approx \left( \sum_{k=1}^n  \int^\infty_0 \big\|\Delta_{te_k}  f\big\|^p_{L^p}\frac{dt}{t^{1+s p}}   \right)^{1/p} ,\quad  \forall  f\in W^{s,p}(\RR^n)\,,    
		\]
		where  $\Delta_{te_k}  f (x)  =  f(x+te_k)-f(x)$, and $e_k$ is the $k$-th vector of the canonical basis in $\RR^n$,    $k=1,\dots,n$.	
		\\
		Thanks to this result, \eqref{5.1} is equivalent to  the following inequality 
		\begin{align}\label{5.1a}
			\sum_{k=1}^n  \int^\infty_0 \big\|\Delta_{te_k}  f\big\|^p_{L^p}\frac{dt}{t^{1+\alpha p}}    \approx  \|f\|^p_{\dot{B}^{\alpha}_{p,p}}  \,.
		\end{align}
		
		Then, we first show that
		\begin{equation} \label{5.1c}
			\sum_{k=1}^n  \int^\infty_0 \big\|\Delta_{te_k}  f\big\|^p_{L^p}\frac{dt}{t^{1+\alpha p}}    \lesssim \|f\|^p_{\dot{B}^{\alpha}_{p,p}}  \,. 
		\end{equation}
		It suffices to prove that
		\begin{equation} \label{5.1b}
			\int^\infty_0 \big\|\Delta_{te_1}  f\big\|^p_{L^p}\frac{dt}{t^{1+\alpha p}}    \lesssim \|f\|^p_{\dot{B}^{\alpha}_{p,p}}  \,. 
		\end{equation}
		Indeed, let $\varphi\in\mathcal{S}(\RR^n)$ be such that ${\rm supp}(\hat{\varphi})\subset \big\{ \frac{1}{2}< |\xi|< 2  \big\}$, $\hat{\varphi}(\xi) \not=0$ in   $\big\{ \frac{1}{4}< |\xi|<1 \big\}$,  $\varphi_j(x)= 2^{-jn}\varphi(2^{-j}x)$ for $j\in\mathbb{Z}$, and  $\displaystyle\sum_{j\in\mathbb{Z}}  \hat{\varphi_j}(\xi) =1$   for $\xi\not=0$. 
		\\
		Next, let us set 
		$$\widehat{\psi}_j(\xi) = \big(e^{it\xi_1}-1\big)  \widehat{\varphi}_j(\xi)\,, \quad \xi=(\xi_1,...,\xi_n) \,.$$ 
		Note that  for any $g\in\mathcal{S}(\RR^n)$ $$\mathcal{F}^{-1}\big\{(e^{it\xi_1}-1)  \widehat{g}\big\} = g(x+te_1)- g(x)  \,,$$
		where $\mathcal{F}^{-1}$ denotes by the inverse Fourier transform.
		\\
		Since ${\rm supp}(\widehat{\varphi}_j)  \cap {\rm supp}(\widehat{\varphi}_l) =\emptyset$ whenever $|l-j|\geq 2$, then we have
		\begin{align}\label{5.2}
			\psi_j * f  &= \psi_j *  \left(\sum_{i\in\mathbb{Z}} \varphi_j \right) * f  = \psi_j *  \big(\varphi_{j-1} + \varphi_j +\varphi_{j+1} \big)  * f \,.
		\end{align}	
		Applying  Young's inequality yields
		\begin{align}\label{5.3}
			\| \psi_j * \varphi_{j} * f \|_{L^p}& \leq  \| \psi_j \|_{L^1} \| \varphi_{j} * f \|_{L^p}  \nonumber
			\\
			&=   \big\| \mathcal{F}^{-1}\big\{ (e^{it\xi_1}-1) \widehat{\varphi}_j(\xi) \big\}\big\|_{L^1} \| \varphi_{j} * f \|_{L^p} \nonumber
			\\
			&=   \| \varphi_j(.+te_1)-\varphi_{j}(.)\|_{L^1} 
			\| \varphi_{j} * f \|_{L^p}  \leq C\| \varphi_{j} * f \|_{L^p}    \,,
		\end{align}	
		where  $C=C_\varphi$ is independent of $j$. 
		\\
		On the other hand, we observe that
		\begin{align*}
			\big|\varphi_j(x+te_1)-\varphi_{j}(x)\big|&=\big| \int^1_0	D\varphi_{j} (x + \tau t e_1)  \cdot te_1  \, d\tau \big| 
			\\
			&\leq   t\int^1_0	\big|D\varphi_{j} (x + \tau t e_1) \big|  \, d\tau = t 2^{-j}  2^{-jn} \int^1_0  	\big|D\varphi \big( 2^{-j}( x + \tau t e_1)\big) \big|  \, d\tau  \,.
		\end{align*}
		Therefore,
		\begin{align}\label{5.5}
			\| \varphi_j(.+te_1)-\varphi_{j}(.)\|_{L^1}  &\leq t 	 2^{-j}  2^{-jn} \int^1_0  	\big\|D\varphi \big( 2^{-j}( x + \tau t e_1)\big) \big\|_{L^1}  \, d\tau \nonumber
			\\
			& =  t 	 2^{-j} \int^1_0 \|D\varphi\|_{L^1}  \, d\tau    = C(\varphi) \, t 2^{-j}  \,.
		\end{align}
		Combining \eqref{5.2},  \eqref{5.3} and \eqref{5.5} yields 
		\begin{equation}\label{5.6}
			\sum_{j\in\mathbb{Z}}	\| \psi_j * f \|_{L^p}	\lesssim   \min\{1,t2^{-j}\}	\sum_{j\in\mathbb{Z}}\| \varphi_{j} * f \|_{L^p}\,.
		\end{equation} 
		Now,  remind that\, $f(x+te_1)-f(x) =  \displaystyle\sum_{j\in\mathbb{Z}}   \psi_j * f(x)$  in $\mathcal{S}^\prime(\RR^n)$. 
		Then, we deduce from  \eqref{5.6} that 
		\begin{align*}
			\int^\infty_0 \int_{\RR^n}   \frac{|f(x+te_1)-f(x)|^p}{t^{1+\alpha p}} \, dx dt  &=  \int_{0}^{\infty} \big\| \sum_{j\in  \mathbb{Z}}  \psi_j * f\big\|_{L^p}^p \frac{dt}{t^{1+\alpha p}}	
			\\
			&\lesssim 	 \sum_{k\in\mathbb{Z}}  \int^{2^k}_{2^{k-1}} \sum_{j\in  \mathbb{Z}}  \min\{1,t^p  2^{-jp}\} \|  \varphi_j * f\|_{L^p}^p \frac{dt}{t^{1+\alpha p}}
			\\
			&\lesssim 	 \sum_{k\in\mathbb{Z}}  2^{-k\alpha p}   \sum_{j\in  \mathbb{Z}}  \min\{1,2^{(k-j)p}\} \|  \varphi_j * f\|_{L^p}^p 
			\\
			&= \sum_{k\in\mathbb{Z}}    \sum_{j\in  \mathbb{Z}}  \min\{1,2^{(k-j)p}\}   2^{-(k-j)\alpha p}  \left[ 2^{-j\alpha p}  \|\varphi_j * f\|_{L^p}^p \right]
			\\
			&= \sum_{k\in\mathbb{Z}}    \sum_{j\in  \mathbb{Z}}  \min\{2^{-(k-j)\alpha p},2^{(k-j)(1-\alpha)p}\}    \left[ 2^{-j\alpha p}  \|\varphi_j * f\|_{L^p}^p \right]
			\\
			&\leq \sum_{k\in\mathbb{Z}}    \sum_{j\in  \mathbb{Z}}  2^{-|k-j|\delta}  \left[ 2^{-j\alpha p}  \|\varphi_j * f\|_{L^p}^p \right],\quad \delta=\min\{\alpha p  ,  (1-\alpha)p\} 
			\\
			&\leq  C_\delta \sum_{k\in\mathbb{Z}} \left[ 2^{-k\alpha p}  \|\varphi_k * f\|_{L^p}^p \right]  = C_\delta \|f\|_{\dot{B}^{\alpha}_{p,p}}^p \,.
		\end{align*}
		Similarly, we also obtain
		\[ 	\int^\infty_0 \int_{\RR^n}   \frac{|f(x+te_k)-f(x)|^p}{t^{1+\alpha p}} \, dx dt \lesssim  \|f\|_{\dot{B}^{\alpha}_{p,p}}^p ,\quad k=2,\dots,n \,.\]
		This yields \eqref{5.1b}. 
		\\
		
		For the converse, let $\{\varphi_j\}_{j\in\mathbb{Z}}$ be the sequence above. By following \cite[page 246]{Grevholm},  we can construct function $\psi\in\mathcal{S}(\RR^n)$ such that  $\hat{\psi}(\xi)  =1$ on $\{1/2 \leq|\xi|\leq 2\}$, and  $\widehat{\psi}=\displaystyle\sum^n_{k=1}  \widehat{h}^{k}$, with $h^{k}\in\mathcal{S}(\RR^n)$ satisfies
		\begin{align}\label{5.8a}
			\sup_{t\in(2^{j-1}, 2^j)}	\big\|\frac{\widehat{h}^k_j(\xi)}{e^{it\xi_k}-1}  \big\|_{L^1}  \leq C  , \quad k=1,\dots,n \,,
		\end{align}
		where $h^k_j(x)  =  2^{-jn}h^k_j(2^{-j}x)$, and  constant $C>0$ is independent of $k, j$.  Actually, we only need $\eqref{5.8a}$ holds for $\big\|\frac{\widehat{h}^k_j(\xi)}{e^{it\xi_k}-1}  \big\|_{\mathcal{M}}$ instead of 
		$\big\|\frac{\widehat{h}^k_j(\xi)}{e^{it\xi_k}-1}  \big\|_{L^1}$, where $\mathcal{M}$ is the space of bounded measures on $\RR^n$, and  $\|\mu\|_{\mathcal{M}}$ is the total variation of $\mu$.
		\\
		Next, from the construction of functions $h^k$, $k=1,\dots,n$, there exists a  universal constant $C_1>0$ such that
		\begin{align*}
			\big\|\mathcal{F}^{-1}\big\{\frac{\widehat{h^k_j}(\xi)}{e^{it\xi_k}-1}\big\} \big\|_{L^1}  \leq C_1 \big\|\frac{\widehat{h}^k_j(\xi)}{e^{it\xi_k}-1}  \big\|_{L^1}\,.
		\end{align*}
		With the last inequality noted, we deduce from \eqref{5.8a} that
		\begin{align}\label{5.8c}
			\sup_{t\in(2^{j-1}, 2^j)}	\big\|\mathcal{F}^{-1}\big\{\frac{\widehat{h^k_j}(\xi)}{e^{it\xi_k}-1}\big\} \big\|_{L^1}  \leq C C_1  , \quad k=1,\dots,n \,.
		\end{align}
		Now, observe  that
		$$h^k_j *f = \mathcal{F}^{-1}\big\{\frac{\widehat{h^k_j}(\xi)}{e^{it\xi_k}-1}\big\} * \Delta_{te_k} f \,.$$ 
		Thus, it follows from the triangle inequality, \eqref{5.8c}, and Young's inequality  that  
		\begin{align}\label{5.9}
			\big\|\psi_j* f\big\|_{L^p} &=   	\big\|\sum_{k=1}^n  h^k_j *f  \big\|_{L^p}   \leq \sum_{k=1}^n \big\| h^k_j *f  \big\|_{L^p}  =   \sum_{k=1}^n\big\| \mathcal{F}^{-1}\big\{\frac{\widehat{h^k_j}(\xi)}{e^{it\xi_k}-1}\big\} * \Delta_{te_k} f  \big\|_{L^p}  \nonumber
			\\
			&\leq \sum_{k=1}^n\big\| \mathcal{F}^{-1}\big\{\frac{\widehat{h^k_j}(\xi)}{e^{it\xi_k}-1}\big\}  \big\|_{L^1}  \big\| \Delta_{te_k} f  \big\|_{L^p}  \nonumber
			\\
			&\lesssim   	 
			\sum_{k=1}^n \big\| \Delta_{te_k} f  \big\|_{L^p}  \,,  \quad \text{for all }  t\in(2^{j-1},2^j)\,.
		\end{align}
		On the other hand, it is clear that $\hat{\psi}(\xi)  \hat{\varphi} (\xi)=\hat{\varphi} (\xi)$ since  ${\rm supp}(\hat{\varphi})\subset \{1/2 \leq|\xi|\leq 2\}$.
		\\
		Hence, we obtain from \eqref{5.9} that
		\begin{align*}
			\| \varphi_{j}  *f \|^p_{L^p}   = 	\| \varphi_{j}  * \psi_j * f \|^p_{L^p} \leq \| \varphi_{j} \|^p_{L^1} \|\psi_j * f\|^p_{L^p}  \lesssim \sum_{k=1}^n \big\|\Delta_{te_k} f  \big\|^p_{L^p}  
		\end{align*}
		for all $t\in (2^{j-1},2^j)$.
		\\
		Thus,
		\begin{align*}
			\sum_{j\in\mathbb{Z}} 2^{-j\alpha p} 	\| \varphi_{j}  *f \|^p_{L^p} &\lesssim 	\sum_{j\in\mathbb{Z}} 2^{-j\alpha p} \sum_{k=1}^n  \fint^{2^j}_{2^{j-1}}   \big\|\Delta_{te_k} f  \big\|^p_{L^p}    \,dt  \lesssim \sum_{k=1}^n \int^\infty_0 \|\Delta_{te_k}  f\|^p_{L^p}   \frac{dt}{t^{1+\alpha p}}  
		\end{align*}
		which yields
		\[ \|f\|^p_{\dot{B}^{\alpha}_{p,p}}  \lesssim  \sum_{k=1}^n \int^\infty_0 \|\Delta_{te_k}  f\|^p_{L^p}   \frac{dt}{t^{1+\alpha p}}  \,. \]
		This completes the proof of Proposition \ref{Pro-cha}.

	\end{proof}

	\bigskip	
	\textbf{Acknowledgement.} The research is funded by University of Economics Ho Chi Minh City, Vietnam.

\end{document}